\documentclass[11pt,reqno]{amsart}
\usepackage{geometry} 
\usepackage{tikz}
\usepackage{tikz-cd}
\usepackage{color}
\usepackage{hyperref}
\hypersetup{
  colorlinks   = true, 
  urlcolor     = blue, 
  linkcolor    = blue, 
  citecolor   = red 
}

\usepackage{amsfonts}
\usepackage{amsmath}
\usepackage{amssymb}
\usepackage{mathtools}
\usepackage{color}
\usepackage{tikz}
\usetikzlibrary{cd}
\usepackage{hyperref}
\usepackage{cancel}

\everymath{\displaystyle}
\numberwithin{equation}{section}

\newtheorem{theorem}{Theorem}[section]
\newtheorem{definition}[theorem]{Definition}
\newtheorem{corollary}[theorem]{Corollary}
\newtheorem{remark}[theorem]{Remark}
\newtheorem{lemma}[theorem]{Lemma}
\newtheorem{proposition}[theorem]{Proposition}

\newtheorem{assumption}[theorem]{Assumption}

\newcommand{\leve}{\big\vert}
\newcommand{\rive}{\big\vert}
\newcommand{\leVe}{\big\Vert}
\newcommand{\riVe}{\big\Vert}

\newcommand{\circg}{\stackrel{\circ}{g}\mkern-1mu \makebox[0ex]{}}    
\newcommand{\circnabla}{ \mkern-2mu \stackrel{\circ}{\nabla} \mkern-5mu \makebox[0ex]{}}    
\newcommand{\circDelta}{\stackrel{\circ}{\Delta} \mkern-5mu \makebox[0ex]{}}     
\newcommand{\circdiv}{\stackrel{\circ}{\div} \mkern-5mu \makebox[0ex]{}}      

\newcommand{\us}{\underline{s}}   
\newcommand{\uC}{\underline{C}}  

\newcommand{\slashg}{g \mkern-8.7mu \slash}  
\newcommand{\dvol}{\mathrm{dvol}}               

\newcommand{\uL}{\underline{L}}

\newcommand{\udelta}{\underline{\delta}}

\renewcommand{\d}{\mathrm{d}}        
\newcommand{\slashd}{\d\mkern-8.5mu /\,}      
\renewcommand{\div}{\mathrm{div}}
\newcommand{\lie}{\mathcal{L}}           

\newcommand{\uh}{\underline{h}}  
\newcommand{\uf}{\underline{f}}    
\newcommand{\uvareps}{\underline{\varepsilon}}    
\newcommand{\dd}[1]{\mathfrak{d}\big\{ #1 \big\}} 
\newcommand{\ufrakd}{\underline{\mathfrak{d}}} 
\newcommand{\bdd}[1]{\delta \big\{ #1 \big\}} 
\newcommand{\er}[1]{\mathbf{e}\big\{#1\big\}} 
\newcommand{\eps}{\epsilon}

\newcommand{\uhl}[1]{\mkern+1mu{}^{#1}\mkern-2.2mu \uh}
\newcommand{\ufl}[1]{\mkern+1mu{}^{#1}\mkern-4mu \uf}
\newcommand{\rel}[1]{\mkern+1mu {}^{#1} \mkern-1mu \mathrm{re}}
\newcommand{\Xl}[1]{\mkern+1mu {}^{#1} \mkern-4mu X}

\newcommand{\Xls}{\mkern+1mu {}^{s} \mkern-4mu X}
\newcommand{\rels}{\mkern+1mu {}^{s} \mkern-1mu \mathrm{re}}

\newcommand{\ufls}{\ufl{s}}
\newcommand{\ucalH}{\underline{\mathcal{H}}}
\newcommand{\dL}{\dot{L}}
\newcommand{\duL}{\dot{\uL}}
\newcommand{\dpartial}{\dot{\partial}}

\newcommand{\dslashd}{\dot{\slashd}}
\newcommand{\db}{\dot{b}}

\newcommand{\dcirc}{\dot{\circ}}
\newcommand{\dcircnabla}{\makebox[2ex]{$\stackrel{\dcirc}{\nabla}$} \makebox[0ex]{}}
\newcommand{\dcircDelta}{\makebox[2ex]{$\stackrel{\dcirc}{\Delta}$} \makebox[0ex]{}}

\newcommand{\Sigmal}[1]{\mkern+1mu {}^{#1} \mkern-1mu \Sigma}

\newcommand{\ud}{\underline{d}}
\newcommand{\ubfd}{\underline{\mathbf{d}}}

\newcommand{\hle}[1]{{}^{#1} \mkern-2mu h}  
\newcommand{\dcircdiv}{\makebox[3.5ex]{$\overset{\dcirc}{\div}$} \makebox[0ex]{}}

\title[Global regular null hypersurfaces in a perturbed Schwarschild exterior]{Global regular null hypersurfaces in a perturbed Schwarzschild black hole exterior}
\author{Pengyu Le}
\address{Yanqi Lake Beijing Institute of Mathematical Sciences and Applications, Beijing, China}
\email{pengyu.le@bimsa.cn}
\date{} 

\begin{document}

\begin{abstract}
The spherically symmetric null hypersurfaces in a Schwarzschild spacetime are smooth away from the singularities and foliate the spacetime. We prove the existence of more general foliations by null hypersurfaces without the spherical symmetry condition. In fact we also relax the spherical symmetry of the ambient spacetime and prove a more general result: in a perturbed Schwarzschild spacetime (not necessary being vacuum), nearly round null hypersurfaces can be extended regularly to the past null infinity, thus there exist many foliations by regular null hypersurfaces in the exterior region of a perturbed Schwarzschild black hole. A significant point of the result is that the ambient spacetime metric is not required to be differentiable in all directions.
\end{abstract}
\maketitle
\tableofcontents

\section{Introduction}\label{section 1}
\noindent
The Schwarzschild black hole spacetime is the static spherically symmetric vacuum solution of the Einstein equations. Its metric reads as follows:
\begin{align}
\nonumber
\d s^2 = - \left( 1-\frac{2m}{r} \right) \d t^2 + \left( 1- \frac{2m}{r} \right)^{-1} \d r^2 + r^2 \left( \d \theta^2 + \sin^2 \theta \d \phi^2 \right).
\end{align}
This metric depends on the parameter $m$. The physical meaning of $m$ is the mass of the spacetime. When $m=0$, it becomes the flat Minkowski metric. At first sight, $r=0$ and $r=2m$ look like the values for which the metric is singular. Only $r=0$ is a true singularity but $r=2m$ is just a coordinate singularity, which can be resolved by coordinate transformations. This fact is revealed by \cite{Ed} \cite{Le} \cite{Fi} \cite{Sy} \cite{Fi} \cite{Kr} \cite{Sz}. In the following, we denote $2m$ by $r_0$.

We employ the following double null coordinate system $\{\us,s,\theta, \phi\}$ of the Schwarzschild black hole spacetime, where the metric takes the form
\begin{align}
\nonumber
g_{S} = \frac{2(s+r_0)}{r}\exp\frac{\us+s+r_0-r}{r_0}\left( \d s\otimes \d \us + \d \us \otimes \d s  \right) + r^2 \left( \d \theta^2 + \sin^2 \theta \d \phi^2 \right),
\end{align}
where
\begin{align}
\nonumber
(r-r_0)\exp\frac{r}{r_0} = s\exp\frac{\us+s+r_0}{r_0}.
\end{align}
In the following, we use $\circg$ to denote the standard metric $\d \theta^2 + \sin^2 \theta \d \phi^2$ on $\mathbb{S}^2$.

We consider a neighbourhood $\mathcal{M}$ of the incoming null hypersurface $\uC_{\us=0}$. 
The general form of a perturbed Schwarzschild metric $g$ is assumed to be
\begin{align}
\nonumber
g=4\Omega^2  \d s  \d \us  + \slashg_{\theta\theta} \big( \d \theta - b^{\theta} \d s\big)^2 + 2\slashg_{\theta\phi} \big( \d \theta - b^{\theta} \d s\big)\big( \d \phi - b^{\phi} \d s\big) + \slashg_{\phi\phi} \big( \d \phi - b^{\phi} \d s\big)^2,
\end{align}
or in a general double null coordinate system $\{\us,s,\theta^1, \theta^2\}$,
\begin{align}
\nonumber
g=2\Omega^2  (\d s \otimes \d \us  + \d \us \otimes \d s) + \slashg_{ij} \big( \d \theta^i - b^i \d s\big) \otimes \big( \d \theta^j - b^j \d s\big).
\end{align}
The spacetime is foliated by the incoming null hypersurfaces $\{\uC_{\us}\}$. Each $\uC_{\us}$ intersects with $C_{s=0}$ at spacelike surface $\Sigma_{s=0,\us}$. By definition $\uC_{\us}$ is extended regularly to the past null infinity.
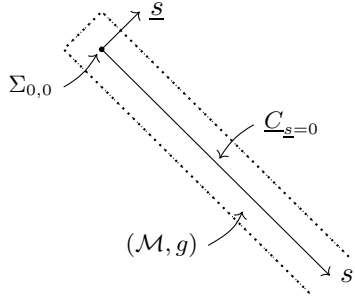
\begin{figure}[h]
\begin{center}
\begin{tikzpicture}
\draw[dotted] (3.25, -2.75) -- (0,0.5) -- (-0.5,0) -- (2.75,-3.25);
\draw[->] (0,0) -- (3,-3) node[right] {$s$};
\draw[->] (0,0) -- (0.5,0.5) node[right] {$\us$};
\draw[dotted,thick] (3.25, -2.75) -- (0,0.5) -- (-0.5,0) -- (2.75,-3.25);
\draw[fill] (0,0)  circle [radius=0.03];
\draw[->] (-0.5,-0.5) node[left] {\footnotesize $\Sigma_{0,0}$} to [out=30,in=-120] (-0.06,-0.06);
\draw[->] (2,-1) node[right] {\footnotesize $\uC_{\us=0}$} to [out=-150,in=60]   (1.55,-1.45);
\draw[->] (1.4,-2.6) node[left] {\footnotesize$(\mathcal{M},g)$}to [out=30,in=-110] (1.85,-2.15);
\end{tikzpicture}
\end{center}
\caption{A perturbed Schwarzschild spacetime $(\mathcal{M},g)$.}
\end{figure}

Suppose that $\Sigma_0$ is an embedded spacelike surface in $C_{s=0}$ and $\ucalH$ is the incoming null hypersurfaces where $\Sigma_0$ is embedded in. In general, we donot expect that $\ucalH$ could be extended regularly to the past null infinity. However, if $\Sigma_0$ is a small perturbation of some surface $\Sigma_{s=0, \us}$, then one could hope that the corresponding null hypersurface $\ucalH$ is a perturbation of the regular null hypersurface $\uC_{\us}$ and also extended regularly to the past null infinity. We confirm the above assertion in this paper.

The strategy to prove that $\ucalH$ is extended regularly to the past null infinity is as follows.
\begin{enumerate}
\item[i.]
Assume that $\ucalH$ is parameterised as the graph of a function $\uh$ in the double null coordinate system $\{s, \us, \vartheta\}$, i.e. $\ucalH$ is parameterised by the following map
\begin{align*}
\nonumber
\varphi_{\uh}: \; (s,\vartheta) \; \mapsto \; (s,\us,\vartheta) = (s, \uh(s,\vartheta), \vartheta) \in \ucalH.
\end{align*}
Suppose that $\ucalH$ intersects with $C_s$ at a spacelike surface $\Sigma_{s}$. Define the function $\ufls$ by
\begin{align*}
\nonumber
\ufls(\vartheta) = \uh(s,\vartheta),
\end{align*}
then $\Sigma_s$ is parametrised by functions $(s, \ufls)$, i.e. $\Sigma_s$ is parametrised by the following map
\begin{align}
\nonumber
\varphi_s:
\;
\vartheta 
\;
\mapsto
\;
(s, \us, \vartheta) = (s, \ufls(\vartheta), \vartheta) \in \Sigma_s.
\end{align}

\item[ii.]
Applying the fact that $\ucalH$ is a null hypersurface, we derive the following first order nonlinear partial differential equations satisfied by $\uh$ and $\ufls$
\begin{align}
\label{eqn 1.7}
\begin{aligned}
\partial_{s} \uh=  -b^i \cdot \uh_i + \Omega^2 \left(\slashg^{-1}\right)^{ij} \cdot \uh_i \cdot \uh_j,
\text{ where }
\uh_i = \partial \uh/\partial \theta^i. 
\\
\partial_{s} \ufls=  -b^i \cdot \ufls_i + \Omega^2 \left(\slashg^{-1}\right)^{ij} \cdot \ufls_i \cdot \ufls_j,
\text{ where } 
\ufls_i = \partial \ufls/\partial \theta^i
\end{aligned}
\end{align}
The converse of the above is also true that if $\uh$ and $\ufls$ solve equations \eqref{eqn 1.7} then $\ucalH$ is null.

\item[iii.] 
Assuming the proper decaying conditions for the metric component $\Omega$, the vector $b^i \partial_i$ and the metric $\slashg$, we can prove the global existence of equation \eqref{eqn 1.7} for small initial data $\ufl{s=0}$. Hence we prove that $\ucalH$ can be extended to the past null infinity regularly.

\end{enumerate}

\begin{figure}[h]
\begin{center}
\begin{tikzpicture}
\draw[dashed] (-1,0)
to [out=70, in=180] (0,0.5)
to [out=0,in=110] (1,0);
\draw (1,0)
to [out=-70,in=0] (0,-0.8) node[below]{\tiny $\Sigma_{0,0}$}
to [out=180,in=-110] (-1,0); 
\draw[dashed] (-1,0) to [out=70,in=-110] (-0.7,0.9);
\draw (-1,0) to [out=-110,in=70] (-2.4,-4);
\draw[dashed] (1,0) to [out=110,in=-70] (0.7,0.9);
\draw[->] (1,0) to [out=-70,in=110] (2.4,-4) node[right] {\small s}; 
\node[above right] at (1.3,-1.3) {\tiny $\uC_{\us=0}$};
\draw[dashed] (-1,0) to [out=-45,in=135] (-0.5,-0.5);
\draw (-1,0) to [out=135, in= -45] (-2,1);
\draw[dashed] (1,0) to [out=-135,in=45] (0.5,-0.5);
\draw[->] (1,0) to [out=45, in= -135] (1.8,0.8) node[right] {\tiny $C_{s=0}$}to [out=45,in=-135] (2.3,1.3) node[right] {\small \us}; 
\draw[blue,thick,->] (1,0) node[right] {\scriptsize  $\ufl{s=0}$} to [out=45, in= -135] (1.65,0.65); 
\draw[blue,thick,->] (-1,0) to [out=135, in= -45] (-1.65,0.65);
\draw[dashed] (-1.5,0.5) to [out=135,in=45] (1.5,0.5);
\draw (1.5,0.5) to [out=-135,in=0] (0,0) node[below]{\scriptsize $\Sigma_0$} to [out=180,in=-45] (-1.5,0.5); 
\draw[dashed] (-1.4,1.3) to [out=-110,in=70] (-1.6,0.7);
\draw (-1.6,0.7) to [out=-110,in=70] (-3,-3);
\draw[dashed] (1.4,1.3) to [out=-70,in=110] (1.6,0.7);
\draw[->] (1.6,0.7) to [out=-70,in=110] (3,-3) node[right]{\small s}; 
\node[right] at (2.1,-0.7) {\scriptsize $\ucalH, \ \us = \uh(s,\vartheta) = \ufls(\vartheta)$};
\draw[dashed] (-1.5,-3.5) -- (-2,-3);
\draw (-2,-3) -- (-3.5,-1.5);
\draw[dashed] (1.5,-3.5) -- (2,-3);
\draw[->] (2,-3) -- (3.5,-1.5); 
\node[right] at (3,-2) {\tiny $C_s$};
\draw[blue,thick,->] (2,-3) to [out=45, in= -135] (2.25,-2.75) node[right] {\scriptsize $\ufls$} to [out=45, in= -135] (2.75,-2.25); 
\draw[blue,thick,->] (-2,-3) to [out=135, in= -45] (-2.75,-2.25);
\draw[dashed] (-2.5,-2.5) to [out=135,in=45] (2.5,-2.5);
\draw (2.6,-2.4) to [out=-135,in=0] (0,-3.4) node[below]{\scriptsize $\Sigma_s$} to [out=180,in=-45] (-2.6,-2.4); 
\end{tikzpicture}
\end{center}
\caption{The null hypersurface $\ucalH$.}
\end{figure}
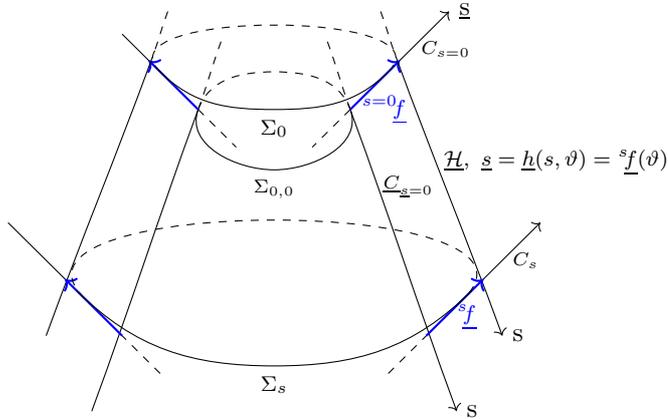

Carrying out the above strategy consists the first part of the paper. The main result of the first part is theorem \ref{thm 3.3}. It should be pointed out that in the perturbed Schwarzschild spacetime considered in theorem \ref{thm 3.3}, it is only assumed that the metric is differentiable tangentially along the outgoing null hypersurfaces $C_s$ but merely continuous in the transversal direction of $\partial_s$, see definition \ref{def 2.2} and remark \ref{rem 2.4}. For example, theorem \ref{thm 3.3} can be applied to spacetimes with impulsive gravitational waves constructed by Luk and Rodnianski in \cite{Lu-Ro}.

The second part of the paper deals with the perturbation of null hypersurfaces in $(\mathcal{M},g)$. The basic problem is as follows. Suppose that two spacelike surfaces $\{\Sigmal{a}_0, a=1,2\}$ are embedded in $C_{s=0}$. Let $\{\ucalH_{a}, a=1,2\}$ be the incoming null hypersurfaces where $\Sigmal{a}$ is embedded in respectively. We want to describe the perturbation between $\ucalH_{a}$. 

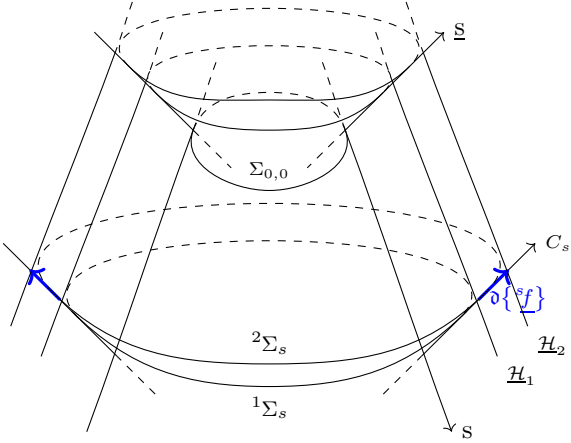
\begin{figure}[h]
\begin{center}
\begin{tikzpicture}
\draw[dashed] (-1,0)
to [out=70, in=180] (0,0.5)
to [out=0,in=110] (1,0);
\draw (1,0)
to [out=-70,in=0] (0,-0.8) node[above]{\tiny $\Sigma_{0,0}$}
to [out=180,in=-110] (-1,0); 
\draw[dashed] (-1,0) to [out=70,in=-110] (-0.7,0.9);
\draw (-1,0) to [out=-110,in=70] (-2.4,-4);
\draw[dashed] (1,0) to [out=110,in=-70] (0.7,0.9);
\draw[->] (1,0) to [out=-70,in=110] (2.4,-4) node[right] {\small s}; 
\draw[dashed] (-1,0) to [out=-45,in=135] (-0.5,-0.5);
\draw (-1,0) to [out=135, in= -45] (-2.3,1.3);
\draw[dashed] (1,0) to [out=-135,in=45] (0.5,-0.5);
\draw[->] (1,0) to [out=45, in= -135] (1.8,0.8) to [out=45,in=-135] (2.3,1.3) node[right] {\small \us}; 
\draw[dashed] (-1.5,-3.5) -- (-2,-3);
\draw (-2,-3) -- (-3.5,-1.5);
\draw[dashed] (1.5,-3.5) -- (2,-3);
\draw[->] (2,-3) -- (3.5,-1.5) node[right] {\tiny $C_s$};
\draw[blue,very thick,->] (2.75,-2.25) node[right] {\scriptsize $\dd{\ufls}$} to [out=45,in=-135] (3.15,-1.85);
\draw[blue,very thick,->] (-2.75,-2.25) to [out=135,in=-45] (-3.15,-1.85);
\draw[dashed] (-1.5,0.5) to [out=135,in=45] (1.5,0.5);
\draw (1.5,0.5) to [out=-135,in=0] (0,0) to [out=180,in=-45] (-1.5,0.5); 
\draw[dashed] (-1.4,1.3) to [out=-110,in=70] (-1.6,0.7);
\draw (-1.6,0.7) to [out=-110,in=70] (-3,-3);
\draw[dashed] (1.4,1.3) to [out=-70,in=110] (1.6,0.7);
\draw (1.6,0.7) to [out=-70,in=110] (3,-3) node[below right] {\tiny $\ucalH_1$}; 
\draw[dashed] (-2.5,-2.5) to [out=135,in=45] (2.5,-2.5);
\draw (2.6,-2.4) to [out=-135,in=0] (0,-3.4) node[below]{\scriptsize $\Sigmal{1}_s$} to [out=180,in=-45] (-2.6,-2.4); 
\draw[dashed] (-1.85,0.85) to [out=135,in=45] (1.85,0.85);
\draw (1.85,0.85) to [out=-135,in=0] (0,0.4)  to [out=180,in=-45] (-1.85,0.85); 
\draw[dashed] (-1.8,1.7) to [out=-110,in=70] (-2,1.1);
\draw (-2,1.1) to [out=-110,in=70] (-3.4,-2.6);
\draw[dashed] (1.8,1.7) to [out=-70,in=110] (2,1.1);
\draw (2,1.1) to [out=-70,in=110] (3.4,-2.6) node[below right] {\tiny $\ucalH_2$}; 
\draw[dashed] (-2.85,-2.15) to [out=135,in=45] (2.85,-2.15);
\draw (2.85,-2.15) to [out=-135,in=0] (0,-3.1) node[above]{\scriptsize $\Sigmal{2}_s$} to [out=180,in=-45] (-2.85,-2.15); 
\end{tikzpicture}
\end{center}
\caption{A perturbation of null hypersurfaces: $\ucalH_1$ and $\ucalH_2$.}
\end{figure}

We still adopt the construction described above. $\uhl{a}$ is the parametrisation function of $\ucalH$. $\Sigmal{a}_s$ is the intersecting spacelike surface of $\ucalH_{a}$ and $C_s$. $(s, \ufl{a,s})$ is the parametrisation functions of $\Sigmal{a}_s$ in the double null coordinate system $\{s,\us,\vartheta\}$. Thus we introduce the perturbation functions
\begin{align}
\nonumber
\dd{\uh}= \uhl{a=2} - \uhl{a=1},
\quad
\dd{\ufls}= \ufl{2,s} - \ufl{1,s} =\dd{\uh}\big|_s .
\end{align}
Then the basic problem of perturbations of null hypersurfaces can be stated as follows: given the null hypersurface $\ucalH_1$, its foliation by $\{\Sigmal{1}_s\}$ and the corresponding parametrisation functions $\uhl{1}$, $\ufl{1,s}$, describe the perturbation functions $\dd{\uh}$, $\dd{\ufls}$ for any $s$ in terms of the initial perturbation function $\dd{\ufl{s=0}}$. It is equivalent to study the perturbation of the solutions of equation \eqref{eqn 1.7} at solutions $\uhl{1}$ and $\ufl{1,s}$. The main results for the perturbation of null hypersurfaces are theorems \ref{thm 4.2} and \ref{thm 4.6}.

Furthermore, we construct an appropriate linearisation for perturbations of null hypersurfaces with the help of equation \eqref{eqn 1.7}. Still adopting the above notations on perturbations of null hypersurfaces, we use $\bdd{\uh}$ and $\bdd{\ufls}$ to denote the linearised perturbation of the parametrisation functions $\uh$ and $\ufls$ respectively, where $\bdd{\uh}\big|_s = \bdd{\ufls}$. We derive the following linear equation of $\bdd{\ufls}$
\begin{align}
\nonumber
&
\overline{\bdd{\ufls}}^{\circg} = \overline{\bdd{\ufl{s=0}}}^{\circg},
\\
\nonumber
&
\partial_s \big( \circDelta \bdd{\ufls} \big) = \Xl{1,s}^i \partial_i \big( \circDelta \bdd{\ufls} \big) - \overline{\Xl{1,s}^i \partial_i \big( \circDelta \bdd{\ufls} \big)}^{\circg},
\end{align}
where $\overline{\makebox[0.25ex]{}\cdot\makebox[0.25ex]{}}^{\circg}$, the overline with a subscript $\circg$, denotes mean value on $\mathbb{S}^2$ with respect to the corresponding area element, $\circDelta$ denotes the Laplacian of $\circg$ on $\mathbb{S}^2$,
the vector field $\Xl{1,s} = \Xl{1,s}^i \partial_i$ is defined by
\begin{align}
\nonumber
\Xl{1,s}^i = - b^i  + 2\Omega^2 \left( \slashg^{-1} \right)^{ij} \ufl{1,s}_j.
\end{align}

We show that the above constructed linearised perturbation $\bdd{\ufls}$ is appropriate by estimating the error $\er{\ufls}$ between the actual perturbation $\dd{\ufls}$ and $\bdd{\ufls}$: $\er{\ufls} = \dd{\ufls} - \bdd{\ufls}$. The main result for the linearised perturbation is theorem \ref{thm 5.2}.

\section{Schwarzschild black hole spacetime and its perturbations}\label{section 2}
\noindent
In this section, we give a short review of the Schwarzschild spacetime in double null coordinate systems and introduce the perturbed Schwarzschild spacetime considered in this paper.

\subsection{The Schwarzschild metric}\label{subsec 2.1}
In the coordinate system $\{t,r,\theta,\phi\}$, the Schwarzschild metric $g_{S}$ takes the following form,
\begin{align}
\nonumber
g_{S}=-\left( 1-\frac{2m}{r} \right) \d t^2 + \left( 1-\frac{2m}{r} \right)^{-1} \d r ^2 + r^{2} \left( \d \theta^2 + \sin^2 \theta \d \phi^2 \right),
\end{align}
where $m$ is a positive number being the total mass of the spacetime.
This form of metric is singular at $r=2m$ and $r=0$. There exist coordinate transformations which resolve the singularity at $r=2m$ and the hypersurface $r=2m$ is an event horizon instead of a singularity. We define that
\begin{align}
\nonumber
r^*= r+2m\log(r-2m),
\end{align}
and the optical functions
\begin{align}
\nonumber
\left\{
\begin{aligned}
&
v=\exp\frac{t+r^*}{4m}=(r-2m)^{\frac{1}{2}} \exp \frac{t+r}{4m}, 
\\
&
w=-\exp \frac{-t+r^*}{4m}= -(r-2m)^{\frac{1}{2}} \exp \frac{-t+r}{4m}.
\end{aligned}
\right.
\end{align}
The coordinate system $\{v, w, \theta, \phi\}$ is called the Kruskal-Szekeres coordinate system and the Schwarzschild metric takes the following form
\begin{align}
\nonumber
g_{S}= -\frac{8m^2}{r} \exp{ \frac{-r}{2m} } \left( \d v \otimes \d w +\d w \otimes \d v \right) + r^2\left(\d \theta^2+\sin^2 \theta \d \phi^2 \right).
\end{align}
The region $\{0 < vw < 2m,0<w\}$ is the black hole region, since it is the set of points which cannot send signals to the future null infinity. The Schwarzschild black hole in the Kruskal-Szekeres coordinate system can be represented by the following picture.

\begin{figure}[h]
\begin{center}
\begin{tikzpicture}
\draw[->] (-3.2,-3.2) -- (3,3) node[anchor=north west]{event horizon $w=0,r=2m$} -- (3.2,3.2) node[right] {$v$};
\draw[->] (3.2,-3.2) -- (-2,2) node[anchor=north east]{$v=0,r=2m$} -- (-3.2,3.2) node[left] {$w$};
\draw[dotted,thick,domain=-2.5:2.5]  plot(\x, {sqrt(3+\x*\x)}) node[above left]{$vw=2m,r=0$};
\draw[dotted,thick,domain=-2.5:2.5]  plot(\x, {-sqrt(3+\x*\x)}) node[below left]{$vw=2m,r=0$};
\draw[domain=0:4]  plot(\x, {-0.3*\x}) node[right]{$t=\text{const}$};
\draw[domain=-1.7:1.7,variable=\y]  plot({sqrt(5+\y*\y)},\y) node[right]{$r=\text{const}>2m$};
\node at (0,1.2) {black hole};
\node at (3.7,0) {exterior region};
\end{tikzpicture}
\caption{Maximal extension of the Schwarzschild black hole spacetime }
\end{center}
\end{figure}
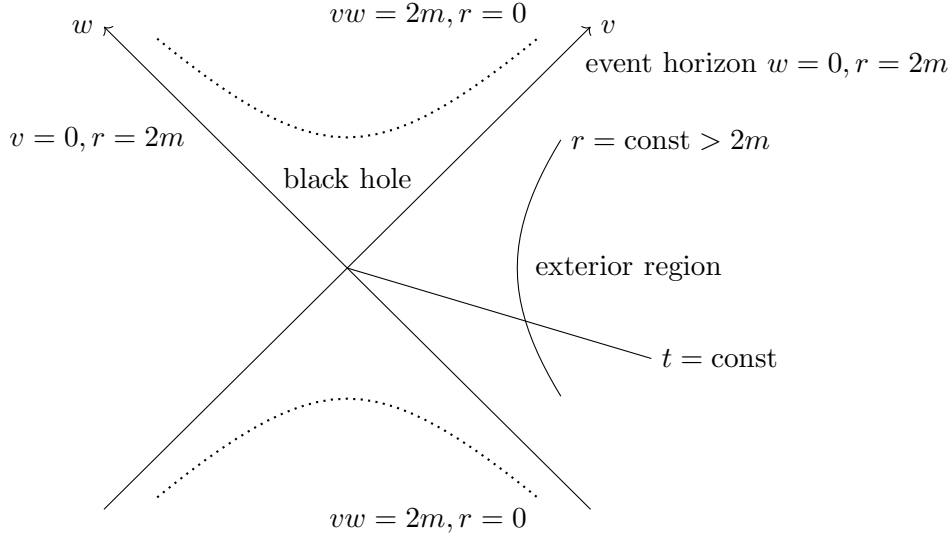

Consider the following coordinate transformation
\begin{align}
\nonumber
\left\{
\begin{aligned}
&
v=v_0\exp\frac{\us}{r_0},
\\
&
w=-\frac{s}{v_0} \exp \frac{s+r_0}{r_0}.
\end{aligned}
\right.
\Leftrightarrow
\left\{
\begin{aligned}
&
(r-2m)^{\frac{1}{2}} \exp \frac{t+r}{4m}=v_0\exp\frac{\us}{r_0},
\\
&
(r-2m)^{\frac{1}{2}} \exp \frac{-t+r}{4m}=\frac{s}{v_0} \exp \frac{s+r_0}{r_0}.
\end{aligned}
\right.
\end{align}
In this coordinate system $\{\us,s,\theta,\phi\}$, the Schwarzschild metric takes the form
\begin{align}
\nonumber
g_{S} = \frac{2(s+r_0)}{r}\exp\frac{\us+s+r_0-r}{r_0}\left( \d s\otimes \d \us + \d \us \otimes \d s  \right) + r^2 \left( \d \theta^2 + \sin^2 \theta \d \phi^2 \right),
\end{align}
where
\begin{align}
\nonumber
(r-r_0)\exp\frac{r}{r_0} = s\exp\frac{\us+s+r_0}{r_0}.
\end{align}
The coordinate systems $\{v,w,\theta,\phi\}$  and $\{\us,s,\theta,\phi\}$ are both double null coordinate systems. We denote the level sets of $s$ by $C_s$, the level sets of $\us$ by $\uC_{\us}$ and use $\Sigma_{\us,s}$ to denote the intersection of $\uC_{\us}$ with $C_s$. 

Associated with the double null coordinate system $\{\us,s,\theta,\phi\}$, we define the corresponding null vectors $L,\uL$ such that $L$ is tangential to the level set $C_s$ and $\uL$ is tangential to the level set $\uC_{\us}$,
\begin{align}
\nonumber
L\us=1,\quad \uL s=1.
\end{align}
We have the following simple formulae for $L,\uL$,
\begin{align}
\nonumber
L=\partial_{\us}, \quad  \uL= \partial_s.
\end{align}
$\{L,\uL\}$ is a null frame of the normal bundle of $\Sigma_{\us,s}$. We define the positive function $\Omega_{S}$ by
\begin{align}
\nonumber
2\Omega_{S}^2=g_{S}\left( L,\uL \right)=\frac{2(s+r_0)}{r} \cdot \frac{r-r_0}{s}.
\end{align}
The null vectors $L',\uL'$ are defined by
\begin{align}
\nonumber
L'=\Omega_{S}^{-2} L, \quad \uL'= \Omega_{S}^{-2} \uL.
\end{align}
We have
\begin{align}
\nonumber
g_{S}(L,\uL') = g_{S}(L', \uL) =2.
\end{align}
Hence $\{L,\uL'\}$ and $\{L',\uL\}$ are conjugate null frame on $\Sigma_{\us,s}$.

\subsection{Perturbations of the Schwarzschild spacetime}\label{subsec 3.3}
We consider a class of perturbed Schwarzschild metric on a neighbourhood of a null hypersurface. We introduce the following definition called the $\kappa$-neighbourhood.
\begin{definition}[$\kappa$-neighbourhood $M_{\kappa}$]
Let $\{\us,s\}$ be the double null foliation of the Schwarzschild spacetime $\left(\mathcal{S},g_{S}\right)$, then the $\kappa$-neighbourhood $M_{\kappa}$ is defined by
\begin{align}
\nonumber
M_{\kappa} = \left\{ p\in \mathcal{S}: s(p)>-\kappa r_0, \leve \us \rive < \kappa r_0 \right\}.
\end{align}
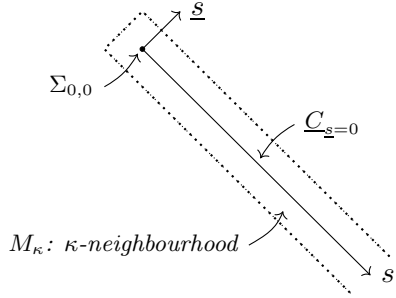
\begin{center}
\begin{tikzpicture}
\draw[dotted] (3.25, -2.75) -- (0,0.5) -- (-0.5,0) -- (2.75,-3.25);
\draw[->] (0,0) -- (3,-3) node[right] {$s$};
\draw[->] (0,0) -- (0.5,0.5) node[right] {$\us$};
\draw[dotted,thick] (3.25, -2.75) -- (0,0.5) -- (-0.5,0) -- (2.75,-3.25);
\draw[fill] (0,0)  circle [radius=0.03];
\draw[->] (-0.5,-0.5) node[left] {\footnotesize $\Sigma_{0,0}$} to [out=30,in=-120] (-0.06,-0.06);
\draw[->] (2,-1) node[right] {\footnotesize $\uC_{\us=0}$} to [out=-150,in=60]   (1.55,-1.45);
\draw[->] (1.4,-2.6) node[left] {\footnotesize$M_{\kappa}$: $\kappa$-neighbourhood }to [out=30,in=-110] (1.85,-2.15);
\end{tikzpicture}
\end{center}
\end{definition}

Let $\{\us,s,\theta,\phi\}$ be the double null coordinates on $M_{\kappa}$ inherited from the Schwarzschild spacetime. A Lorentzian metric on $M_{\kappa}$ has the following form
\begin{align}
\nonumber
g=4\Omega^2  \d s  \d \us  + \slashg_{\theta\theta} \big( \d \theta - b^{\theta} \d s\big)^2 + 2\slashg_{\theta\phi} \big( \d \theta - b^{\theta} \d s\big)\big( \d \phi - b^{\phi} \d s\big) + \slashg_{\phi\phi} \big( \d \phi - b^{\phi} \d s\big)^2,
\end{align}
where all the metric coefficients are functions of the coordinates. We define 
\begin{align}
\nonumber
\slashg=\slashg_{\theta\theta} \d \theta^2+ 2\slashg_{\theta\phi} \d \theta\d \phi+ \slashg_{\phi\phi} \d \phi^2,
\quad
\vec{b}=b^{\theta}\partial_{\theta} + b^{\phi} \partial_{\phi}.
\end{align}
$\slashg$ is the restriction of $g$ on $\Sigma_{\us,s}$ and $\vec{b}$ is a vector field on $\Sigma_{\us,s}$. For the Schwarzschild metric, $\slashg_{S}= r^2 \left( \d \theta^2 + \sin^2\theta \d \phi^2\right)= r^2 \circg$ where $\circg$ is the round metric of radius $1$ on the sphere. We will use $\circnabla$ to denote the covariant derivative of $\circg$ on $\Sigma_{\us,s}$.

Let $\{\theta^1,\theta^2\}$ be any coordinate system of $\mathbb{S}^2$. We obtain a double null coordinate system $\{\us,s,\theta^1,\theta^2\}$, and the metric $g$ takes the following form in this coordinate system
\begin{align}
g=2\Omega^2  ( \d s  \otimes \d \us  +\d \us \otimes \d s ) + \slashg_{ab} \big( \d \theta^a - b^{a} \d s\big) \otimes \big( \d \theta^b - b^{b} \d s\big),
\end{align}
where $a,b$ are indices in $\{1,2\}$. 

By the construction, $\left( M_{\kappa},g\right)$ is a Lorentzian manifold and $\{\us,s,\theta^1,\theta^2\}$ is a double null coordinate system of this manifold. Let $L$ be the null vector field tangent to $C_s$, $\uL$ be the null vector field tangent to $\uC_{\us}$ and
\begin{align}
\nonumber
L\us=1, \quad \uL s=1.
\end{align}
We have
\begin{align}
\nonumber
L=\partial_{\us}, \quad  \uL= \partial_s + \vec{b}.
\end{align}
We can define $L'$ and $\uL'$ by
\begin{align}
\nonumber
L'=\Omega^{-2} L, \quad \uL'=\Omega^{-2} \uL.
\end{align}

In the following, we introduce a class of Lorentzian metrics on $M_{\kappa}$, which are close to the Schwarzschild metric.
\begin{definition}[$\eps$-close Schwarzschild metric $g_{\eps}$]\label{def 2.2}
Let $g_{\eps}$ be a continuous Lorentzian metric on $M_{\kappa}$ that in coordinate system $\{\us,s,\theta^1,\theta^2\}$
\begin{align}
\nonumber
g_{\eps}=2\Omega^2  \left( \d s  \otimes \d \us  +\d \us \otimes \d s \right) + \slashg_{ab} \big( \d \theta^a - b^{a} \d s\big) \otimes \big( \d \theta^b - b^{b} \d s\big).
\end{align}
We define the area radius function $r(\us,s)$ by
\begin{align}
\nonumber
4\pi r^2(\us,s) =\int_{\Sigma_{\us,s}} 1 \cdot \dvol_{\slashg}.
\end{align}
We call $g_{\eps}$ $\eps$-close to the Schwarzschild metric on $M_{\kappa}$, if the following assumptions on the metric $g_{\eps}$ hold:
\begin{align}
\nonumber
&  1-\eps < \frac{r}{r_{S}} < 1+ \eps, 
\\
&
\begin{aligned}
&
\leve \log \Omega - \log \Omega_{S} \rive < \frac{\eps r_0}{r},
&& \leve \circnabla^k \left( \log \Omega - \log \Omega_{S} \right) \rive < \frac{\eps r_0}{r},
\nonumber
\\
\nonumber
&  \leve \circnabla^k \partial_{\us}^m \left( \log \Omega - \log \Omega_{S} \right) \rive < \frac{\eps }{r r_0^{m-1}},
\end{aligned}
\\
&
\begin{aligned}
& \leve \vec{b} \rive_{\circg} < \frac{\eps r_0 \vert \us\vert}{r^3},  
&& \leve \circnabla^{k} \vec{b} \rive_{\circg} < \frac{\eps r_0 \vert\us\vert}{r^3}, 
&& \leve \circnabla^{k} \partial_{\us}^{m} \vec{b} \rive_{\circg} < \frac{\eps }{r^{3}r_0^{m-2}},
\nonumber
\\
&  \leve \slashg - \slashg_{S} \rive_{\circg} < \eps r^2, 
&&  \leve \circnabla^{k} \left( \slashg - \slashg_{S} \right) \rive \leq \eps r^2, 
&& \leve \circnabla^k \partial_{\us}^m \left( \slashg - \slashg_{S} \right) \rive_{\circg} < \frac{\eps r^2}{r_0^{m}}.
\nonumber
\end{aligned}
\end{align}
\end{definition}

\begin{remark}
The decay assumptions in the above definition of $g_{\eps}$ are motived from the results in the nonlinear stability of Minkowski spacetime \cite{CK} by Christodoulou and Klainerman. It is worth to point out that the result \cite{CK} was generalised by Bieri \cite{Bi1}-\cite{Bi3} for the Einstein vacuum equations by relaxing the decay of the initial data and the resulting spacetimes in \cite{Bi1}-\cite{Bi3} have much weaker asymptotic decay rates than those in \cite{CK}.
\end{remark}

\begin{remark}\label{rem 2.4}
Note that there is no assumption on the derivatives of the metric components with respect to $\partial_s$. Thus we only assume that the metric is differentiable in tangential directions on each outgoing null hypersurface $C_s$, but merely continuous in the transversal $\partial_s$ direction.
\end{remark}

In the rest of this paper, we will work in the spacetime $(M_{\kappa},g_{\eps})$. To simplify the notations, we will use $(M,g)$ to denote $(M_{\kappa},g_{\eps})$.

\section{Parametrisation of incoming null hypersurfaces}\label{section 3}
\noindent
In this section, we introduce a method to parametrise incoming null hypersurfaces and derive the equation \eqref{eqn 1.7} satisfied by parametrisation functions. Then we estimate the parametrisation functions for a class of null hypersurfaces and prove that they can be extended regularly to the past null infinity.

\subsection{The equation of parametrisation functions}\label{subsec 3.1}
Let $\ucalH$ be an incoming null hypersurface embedded in $(M,g)$ which is regular. Assume that $\ucalH$ is parametrised by a function $\uh$ as its graph of $\us$ over $(s,\vartheta)$ domain $(-\kappa r_0, + \infty) \times \mathbb{S}^2$ in the double null coordinate system, i.e.
\begin{align*}
\varphi_{\uh}: \; (s,\vartheta) \; \mapsto \; (s,\us,\vartheta) = (s, \uh(s,\vartheta), \vartheta) \in \ucalH.
\end{align*}
We define $\Sigma_s$ as the intersection of $\ucalH$ with $C_s$, i.e.
\begin{align}
\nonumber
\Sigma_s= \ucalH \cap C_s.
\end{align}
In particular, $\Sigma_0$ is the intersection of $\ucalH$ with $C_{s=0}$. $\ucalH$ is determined by $\Sigma_0$. Define the function $\ufls$ by $\ufls(\vartheta) = \uh(s,\vartheta)$, then $\Sigma_s$ is parametrised by $\left(s, \ufls \right)$, i.e.
\begin{align*}
\varphi_s:
\;
\vartheta 
\;
\mapsto
\;
(s, \us, \vartheta) = (s, \ufls(\vartheta), \vartheta) \in \Sigma_s.
\end{align*}

Geometrically, $\{ \Sigma_s \}$ is a foliation of $\ucalH$. Associated with the foliation $\{\Sigma_s\}$, we define the tangential null vector fields $\duL$ on $\ucalH$ by the condition
\begin{align}
\nonumber
\duL s =1,
\end{align}
and define the conjugate null normal vector field $\dL'$. As the convention adopted in this paper, the dotted symbol means that the corresponding object lives on $\Sigma_s$.

We restrict coordinate functions $\{ s, \theta^1, \theta^2 \}$ on $\ucalH$, then obtain a coordinate system of $\ucalH$. 
$\left\{ \dpartial_s, \dpartial_1, \dpartial_2 \right\}$ is the coordinate frame of $\{ s, \theta^1, \theta^2 \}$ on $\ucalH$. They are given by
\begin{align}
\dpartial_s = \partial_s + (\partial_s \uh) \partial_{\us}, 
\quad 
\dpartial_i = \partial_i + (\partial_i \uh) \partial_{\us}.
\label{eqn 3.4}
\end{align}
Since $\uh$ is independent of $\us$, we have $\dpartial_s \uh= \partial_s \uh$, $\dpartial_i \uh = \partial_i \uh$. Similarly $\dpartial_i \ufls= \partial_i \ufls$.

Note that $\dpartial_s$ may not be null, which implies that $\duL$ and $\dpartial_s$ may not coincide, thus we can define a $\Sigma_s$-tangential vector field $\vec{\db}$ by
\begin{align}
\duL = \dpartial_s + \vec{\db}, \quad \vec{\db}=\db^1 \dpartial_1 + \db^2 \dpartial_2.
\label{eqn 3.5}
\end{align}

Applying the calculations in \cite{L2}, we have that $\duL$ is
\begin{align}
\duL= \uL + \uvareps L + \uvareps^i \partial_i = \partial_s + \uvareps \partial_{\us} + ( b^i + \uvareps^i )\partial_i,
\label{eqn 3.6}
\end{align}
where $\uvareps$ and $\uvareps^i$ are given by
\begin{align}
\uvareps = - \Omega^2 \left( \slashg^{-1} \right)^{ij} \uh_i \uh_j = - \Omega^2 \leve \dslashd \uh \rive_{\slashg}^2,  \quad \uvareps^k= -2\Omega^2 \uh_i \left( \slashg^{-1} \right)^{ik}.
\label{eqn 3.7}
\end{align}

We substitute \eqref{eqn 3.4} into \eqref{eqn 3.5}
\begin{align}
\duL=\dpartial_s + \db^i \dpartial_i =
\partial_s + \left( \partial_s \uh + \db^i \dpartial_i \uh \right) \partial_{\us} + \db^i\partial_i,
\label{eqn 3.8}
\end{align}
then from equations \eqref{eqn 3.6} \eqref{eqn 3.8}, we obtain
\begin{align}
&
\db^i = b^i + \uvareps^i,
\label{eqn 3.9}
\\
\nonumber
&
\dpartial_{s} \uh + (b^i + \uvareps^i) \dpartial_i \uh= \uvareps,
\end{align}
and substituting from equation \eqref{eqn 3.7}, we get
\begin{align}
\nonumber
\dpartial_{s} \uh=  -b^i  \dpartial_i \uh + \Omega^2 \left(\slashg^{-1}\right)^{ij} \dpartial_i \uh \dpartial_j \uh.
\end{align}
and fully written out, the equation reads
\begin{align*}
\partial_{s} \uh=  -b^i|_{\us=\uh(s,\vartheta)}  \partial_i \uh + \big( \Omega^2 \left(\slashg^{-1}\right)^{ij} \big)\big|_{\us=\uh(s,\vartheta)} \partial_i \uh \partial_j \uh.
\end{align*}
Introduce the function $\phi=\uh(s,\vartheta)-\us$, then the above equation is the eikonal equation on the $0$-level set of $\phi$.

Thus by the above derivation, we prove the following proposition.
\begin{proposition}\label{pro 3.1}
Let $\ucalH$ be an incoming null hypersurface in $(\mathcal{M},g)$. Suppose that $\uh$ is the parametrisation function $\ucalH$, then $\uh$ satisfies the first order partial differential equation
\begin{align*}
\dpartial_{s} \uh=  -b^i  \dpartial_i \uh + \Omega^2 \left(\slashg^{-1}\right)^{ij} \dpartial_i \uh \dpartial_j \uh.
\end{align*}
Suppose that $\ufls$ is the parametrisation function of the foliation $\{\Sigma_s\}$ on $\ucalH$, then $\ufls$ satisfies the first order partial differential equation
\begin{align}
\dpartial_{s} \ufls=  -b^i  \dpartial_i \ufls + \Omega^2 \left(\slashg^{-1}\right)^{ij} \dpartial_i \ufls \dpartial_j \ufls.
\label{eqn 3.12}
\end{align}
\end{proposition}

The converse of proposition \ref{pro 3.1} is also true.
\begin{proposition}\label{pro 3.2}
Let $\{\Sigma_s\}$ be a family of spacelike surfaces parametrised by $(s, \ufls)$. Suppose that $\ufls$
satisfies equation \eqref{eqn 3.12}, then the hypersurface $\ucalH$ foliated by $\Sigma_s$ is an incoming null hypersurface.
\end{proposition}

Thus by propositions \ref{pro 3.1} and \ref{pro 3.2}, we know that in order to show that an incoming null hypersurface $\ucalH$ can be extend to the past null infinity, it is sufficient to prove that the solution of equation \eqref{eqn 3.12} exists for all $s>s_0$ for some $s_0$. We shall show in section \ref{subsec 3.2}, given the initial data $\ufl{s=0}$ sufficiently small, there exists a global solution of equation \eqref{eqn 3.12}. We emphasise that while equation \eqref{eqn 3.12} is quadratic in the first derivatives of $\uh$, it is fully nonlinear in $\ufls$ itself.

\subsection{Estimates for parametrisation functions}\label{subsec 3.2}
We adopt the notations in the previous section. $\Sigma_0$ is parametrised by $(0,\ufl{s=0})$ and embedded in the incoming null hypersurface $\ucalH$. By definition, if $\ufl{s=0}\equiv \us$, then $\Sigma_0$ is simply the surface $\Sigma_{0,\us}$, and $\ucalH$ is the null hypersurface $\uC_{\us}$ in the double null foliation. By definition \ref{def 2.2}, $\uC_{\us}$ is extended to the past null infinity.

For a spacelike surface $\Sigma_0$ with an arbitrary parametrisation $(0,\ufl{s=0})$, the incoming null hypersurface $\ucalH$ cannot be regularly extended to the past null infinity in general, but for a perturbed spacelike surface $\Sigma_0$ of $\Sigma_{0,\us}$, it is possible to regularly extend $\ucalH$ to the past null infinity. We have the following theorem.
\begin{theorem}\label{thm 3.3}
Let $\Sigma_0$ be a spacelike surface parametrised by $(0,\ufl{s=0})$ in $(M,g)$. Assume that it is embedded in the incoming null hypersurface $\ucalH$. Suppose that the parametrisation function $\ufl{s=0}$ satisfies the following inequalities
\begin{align}
\nonumber
\leVe \dslashd \ufl{s=0} \riVe^{n+1,p}_{\circg} \leq \udelta_o r_0,
\quad
\leve \overline{\ufl{s=0}}^{\circg} \rive \leq \udelta_m r_0,
\end{align}
where $n\geq 1,p>2$ or $n\geq 2, p>1$.

There exists a small positive constant $\delta$ depending on $n,p$, that if $\eps, \udelta_o, \udelta_m$ are suitably bounded such that $\eps, \udelta_o, \eps \udelta_m \leq \delta$, then the incoming null hypersurface $\ucalH$ can be extended regularly to the past null infinity, as long as $\ucalH$ does not reach the boundary null hypersurfaces $\uC_{\us=\kappa r_0}$ or $\uC_{\us=-\kappa r_0}$. Here $\eps$ is the one in definition \ref{def 2.2} of $\eps$-close Schwarzschild metric $g$.

Moreover, let $\{\Sigma_{s}\}$ be the foliation on $\ucalH$ and be parametrised by $(s, \ufls)$. Let $(s_{\alpha},s_{\beta})$ be an interval containing $s=0$. Suppose that $\leve \ufls \rive <\kappa r_0$ for all $s \in (s_{\alpha},s_{\beta})$, then there exist constants $c_o, c_{m,m}, c_{m,o}$ depending on $n,p$ such that $\ufls$ satisfies the following estimates
\begin{align}
\label{eqn 3.14}
&
\leVe \dslashd \ufls \riVe^{n+1,p} 
\leq
c_o \leVe \dslashd \ufl{s=0} \riVe^{n+1,p} 
\leq
c_o \udelta_o r_0,
\\
&
\leve \overline{\ufls}^{\circg} \rive 
\leq
c_{m,m} \leve \overline{\ufl{s=0}}^{\circg} \rive + \frac{c_{m,o}}{r_0} \left( \leVe \dslashd \ufl{s=0} \riVe^{n+1,p} \right)^2 
\leq 
( c_{m,m} \udelta_m + c_{m,o} \udelta_o^2 ) r_0,
\label{eqn 3.15}
\end{align}
for all $s\in (s_{\alpha}, s_{\beta})$.
\end{theorem}

Note that theorem \ref{thm 3.3} requires the assumption that $\ucalH$ does not reach the boundary null hypersurfaces. After the theorem is proved, we shall give a condition on $\udelta_m,\udelta_o$ such that $\ucalH$ never reaches the boundary null hypersurfaces in corollary \ref{coro 3.6} of the theorem.

In order to prove theorem \ref{thm 3.3}, it is sufficient to show that there exist constants $\delta$ and $c_o, c_{m,m}, c_{m,o}$ such that equations \eqref{eqn 3.14} \eqref{eqn 3.15} hold.

We take equation \eqref{eqn 3.12} to estimate $\ufls$. We differentiate the equation to obtain the equation of $\dcircDelta \ufls$.
\begin{align}
\label{eqn 3.16}
\dpartial_s \dcircDelta \ufls
=&
-b^i (\dcircDelta \ufls)_i 
+ 2\Omega^2 {\slashg^{-1}}^{ij} \ufls_j (\dcircDelta\ufls)_i
- b^i \, \ufls_i 
+ 2\Omega^2 {\slashg^{-1}}^{ij}\ \ufls_i \ufls_j
\nonumber
\\
&
- (\circDelta b)^i \ \ufls_i 
- 2 (\circnabla^k b^i) \dcircnabla^2_{ki} \ufls
- \partial_{\us} ( \circnabla^k b^i )\ \ufls_k \ufls_i 
- \partial_{\us} b^i\ \ufls_i \dcircDelta \ufls 
\nonumber
\\
&
- 2 \partial_{\us}b^i\ \ufls^k \dcircnabla^2_{ki} \ufls 
- \partial_{\us}^2 b^i\ \ufls_i \ufls_k \ufls^k
\nonumber
\\
&
+ \circDelta( \Omega^2 \slashg^{-1})^{ij}\ \ufls_i \ufls_j 
+ 4 \circnabla^k( \Omega^2 {\slashg^{-1}})^{ij}\ \ufls_i \dcircnabla^2_{kj} \ufls
\nonumber
\\
&
+ 2\Omega^2 {\slashg^{-1}}^{ij} \dcircnabla_k \ufls_i \dcircnabla^k \ufls_j 
+2 \ufls_k \partial_{\us}\circnabla^k ( \Omega ^2 \slashg^{-1} )^{ij} \ufls_i \ufls_j 
\nonumber
\\
&
+ 4 \partial_{\us}(\Omega^2 \slashg^{-1})^{ij} \ufls_i \ufls^k \dcircnabla^2_{kj} \ufls
+ \partial_{\us} ( \Omega^2 \slashg^{-1})^{ij} \ufls_i \ufls_j  \dcircDelta \ufls
\nonumber
\\
&
+ \partial_{\us}^2 ( \Omega^2 \slashg^{-1})^{ij} \ufls_i \ufls_j \ufls^k \ufls_k.
\end{align}
$\dcircDelta, \dcircnabla$ are the Laplacian operator and covariant derivative on $(\Sigma_s, \circg)$, while $\circDelta, \circnabla$ are the Laplacian operator and covariant derivative on $(\Sigma_{s,\us}, \circg)$,
\begin{align}
\nonumber
\dcircDelta = \frac{1}{ \sqrt{\det (\circg_{ab}) }} \dpartial_i \left( {\circg^{-1}}^{ij} \sqrt{\det (\circg_{ab})} \dpartial_i \right),
\;
\circDelta = \frac{1}{ \sqrt{\det (\circg_{ab}) }} \partial_i \left( {\circg^{-1}}^{ij} \sqrt{\det (\circg_{ab})} \partial_i \right).
\end{align}

We can rephrase equation \eqref{eqn 3.16} as follows.
\begin{align}
\label{eqn 3.18}
\dpartial_s \dcircDelta \ufls
=&
-b^i (\dcircDelta \ufls)_i + 2\Omega^2 {\slashg^{-1}}^{ij} \ufls_j (\dcircDelta\ufls)_i
\nonumber
\\
&
- b^i \, \ufls_i
+ 
\sum_{\substack{m+n=3\\ 2\geq m'+n \geq 1}} \circnabla^n \partial_{\us}^{m'} b  \ast P_{m,m'+1}( \dcircnabla, \ufls) 
\nonumber
\\
&
+ 
\sum_{\substack{m+n=4 \\ 2\geq m'+n \geq 1}} P_{m,m'+2}( \dcircnabla, \ufls ) \ast \circnabla^{n} \partial_{\us}^{m'} ( \Omega^2 \slashg^{-1} ) 
\nonumber
\\
& 
+ \sum_{\substack{m_1+m_2 \leq 4 \\ m_1 \leq 2, m_2 \leq 2}} P_{m_1,1} ( \dcircnabla, \ufls) \ast P_{m_2,1}(\circnabla, \ufls) \ast (\Omega^2 \slashg^{-1} ).
\end{align}
Here $P_{m,m'}(\dcircnabla, \ufls)$ means polynomials of $\dcircnabla^l \ufls$ where $m$ is the sum of the order of covariant derivatives and $m'$ is the degree of the polynomial, i.e.
\begin{equation}
\nonumber
P_{m,m'}( \dcircnabla, \ufls) = \sum_{\substack{ k_1 + \cdots + k_{m'} =m \\ k_1, \cdots, k_{m'}\geq 1}} c_{k_1\cdots k_{m'}} \Big( \dcircnabla^{k_1} \ufls \Big) \cdots \Big( \dcircnabla^{k_{m'}} \ufls \Big).
\end{equation}
For example $\ufls_i \dcircnabla_{jk} \ufls$ is a polynomial of type $P_{3,2}(\dcircnabla,\ufls)$. $\ast$ means the contraction of tensor fields with respect to the metric $\circg$. For example we write $- 2 \partial_{\us}b^i\ \ufls_k \dcircnabla^k \dcircnabla_i \ufls $ as $\partial_{\us} b \ast \dcircnabla \ufls \ast \dcircnabla^2 \ufls$.

In order to estimate $\dcircDelta \ufls$ by equations \eqref{eqn 3.16} \eqref{eqn 3.18}, we emphasis the structure of these equations by rewriting them as follows.
\begin{align}
\label{eqn 3.20}
&
\dpartial_s \dcircDelta \ufls = \Xls^i \dpartial_i \big( \dcircDelta \ufls \big) + \rels.
\end{align}
Here $\Xls$ is the vector field
\begin{align}
\Xls= \Xls^i \dpartial_i, 
\quad 
\Xls^i = -b^i + 2\Omega^2 (\slashg^{-1})^{ij} \ufls_j,
\label{eqn 3.21}
\end{align}
and $\rels$ is the rest terms remained in equations \eqref{eqn 3.16} \eqref{eqn 3.18}. Note that $\Xls^i=-\db^i$ by equations \eqref{eqn 3.7} \eqref{eqn 3.9}, 
then by equation \eqref{eqn 3.5}, 
we have $\duL = \dpartial_s - \Xls^i \dpartial_i$. Thus equation \eqref{eqn 3.20}
is simply
\begin{align}
\nonumber
\duL \dcircDelta \ufls = \rels.
\end{align}

In the structure of equation \eqref{eqn 3.20}, it is crucial that $\Xls$ involves only the first order derivative of $\ufls$, and the highest order derivatives of $\ufls$ in $\rels$ is of second order.

We shall use Gronwall's inequality, together with the elliptic estimates on the sphere, to control the Sobolev norm $\leVe \dslashd \ufls \riVe^{n+1,p}$ for the derivatives of $\ufls$. The nonlinear term $ 2\Omega^2 \left(\slashg^{-1}\right)^{ij} \dcircnabla_k \ufls_i \dcircnabla^k \ufls_j$ in $\rels$ is worth to be payed more attentions, since it is a quadratic term of the highest order derivatives. We need the Sobolev norm of sufficiently high order derivatives of $\ufls$ in order to control the $L^{\infty}$ norm of $\dcircnabla^2 \ufls$, such that Gronwall's inequality gives the desired estimate. This is the reason that we choose the ranges of $n,p$ as stated in the assumption of theorem \ref{thm 3.3}.

We explain more on the structure of $\rels$. $\rels$ is a quadratic nonlinear term. Formally, we assume that $\dcircnabla^k \ufls$ and $\circnabla^k \partial_{\us}^l b$ are of size $\delta$, then $\rels$ is of size $\delta^2$. This leads to the assumption on the bounds of $\eps, \udelta_o, \udelta_m$ in theorem \ref{thm 3.3}. 
Note that the smallness condition on $\udelta_m$ is $\eps \udelta_m<\delta$ instead of $\udelta_m<\delta$. This condition is necessary to bound $b, \circnabla^k b$ because they are bounded by $\frac{\eps r_0 |\us|}{r^3}$ in definition \ref{def 2.2}.

We prove theorem \ref{thm 3.3} by bootstrap arguments. We introduce the following bootstrap assumption first.
\begin{assumption}[Bootstrap assumption of $\ufls$ for theorem \ref{thm 3.3}]\label{assum 3.4}
Estimates \eqref{eqn 3.14} \eqref{eqn 3.15} of $\ufls$ hold for all $s\in [0,s_a] \subseteq (s_{\alpha},s_{\beta})$ or $s\in[s_a,0]\subseteq (s_{\alpha},s_{\beta})$.
\end{assumption}

The basic idea of the bootstrap argument is to use the bootstrap assumption to improve the inequalities in it such that the bootstrap assumption can be extended further, and eventually hold for the whole range. We shall use the bootstrap assumption to estimate the vector field $\Xls$ and the nonlinear terms $\rels$, then apply Gronwall's inequality to improve the estimates of $\ufls$ in the bootstrap assumption, such that the bootstrap argument gives the proof of theorem \ref{thm 3.3}.

\begin{proof}[Proof of theorem \ref{thm 3.3}]
We use $c(n,p)$ to denote constants depending only on $n,p$, which not necessary to be same in the proof. We first assume that $\delta\leq \frac{1}{2}$, then $\eps, \udelta_o, \eps \udelta_m \leq \frac{1}{2}$. Moreover, we assume that $\delta$ is sufficiently small such that
\begin{align}
(c_o+c_{m,m} +c_{m,o}) \delta \leq 1.
\nonumber
\end{align}
Assume the bootstrap assumption \ref{assum 3.4} holds, we shall prove the following estimates for $\Xls$ and $\rels$,
\begin{align}
\leVe \Xls \riVe^{n+1,p}
\leq
&
\frac{c(n,p) \eps r_0}{(r_0+s)^3} \left( \leve \overline{\ufls}^{\circg} \rive + \leVe \dslashd \ufls \riVe^{n,p} \right)
+
\frac{c(n,p)}{(r_0+s)^2} \leVe \dslashd \ufls \riVe^{n+1,p}
\nonumber
\\
\leq
&
\frac{c(n,p)( c_{m,m} \eps \udelta_m + c_{m,o} \eps \udelta_o^2 + c_o \eps \udelta_o )r_0^2}{(r_0+s)^3}
+
\frac{c(n,p) c_o \udelta_o r_0}{(r_0+s)^2},
\label{eqn 3.23}
\\
\leVe \rels \riVe^{n,p}
\leq
&
\frac{c(n,p) \eps r_0}{(r_0+s)^3} \left( \leve \overline{\ufls}^{\circg} \rive + \leVe \dslashd \ufls \riVe^{n,p} \right) \leVe \dslashd \ufls \riVe^{n+1,p}
+
\frac{c(n,p)}{(r_0+s)^2} \left( \leVe \dslashd \ufls \riVe^{n+1,p} \right)^2
\nonumber
\\
\leq
&
\frac{c(n,p)( c_{m,m} \eps \udelta_m + c_{m,o} \eps \udelta_o^2 + c_o \eps \udelta_o )c_o \udelta_o r_0^3}{(r_0+s)^3}
+
\frac{c(n,p) c_o^2 \udelta_o^2 r_0^2}{(r_0+s)^2}
\label{eqn 3.24}
\end{align}
First, by bootstrap assumption \ref{assum 3.4} and the Sobolev embedding on the sphere, we have the $L^{\infty}$ estimate of $\ufls$
\begin{align}
&
\leve \ufls \rive 
\leq
c(n,p) (c_{m,m} \udelta_m + c_{m,o} \udelta_o^2 + c_o \udelta_o ) r_0,
\nonumber
\\
&
\leve \dcircnabla^{m} \ufls \rive 
\leq
c(n,p) c_o \udelta_o r_0,
\quad
m\leq n \text{ if } p>2, \text{ or } m \leq n-1, \text{ if } 2\geq p > 1.
\nonumber
\end{align}
Therefore by the $L^{\infty}$ bounds of $b^i$, $\Omega$, $\slashg$ and their derivatives in definition \ref{def 2.2}, we have
\begin{align}
&
\leVe b^i \dpartial_i \riVe^{n+1,p}
\leq
\frac{c(n,p) \eps ( c_{m,m} \udelta_m + c_{m,o} \udelta_o^2 + c_o \udelta_o )r_0^2}{(r_0+s)^3},
\nonumber
\\
&
\leVe 2\Omega^2 (\slashg^{-1})^{ij} \ufls_j \dpartial_i \riVe^{n+1,p}
\leq
\frac{c(n,p) c_o \udelta_o r_0}{(r_0+s)^2},
\nonumber
\end{align}
hence estimate \eqref{eqn 3.23} follows.

Similarly in $\rels$, we have
\begin{align}
&
\sum_{\substack{m+n=3 \\ 2\geq m'+n \geq 1}}
\leVe \circnabla^n \partial_{\us}^{m'} b  \ast P_{m,m'+1}( \dcircnabla, \ufls)  \riVe^{n,p}
\leq
\frac{c(n,p) \eps ( c_{m,m} \udelta_m + c_{m,o} \udelta_o^2 + c_o \udelta_o )c_o \udelta_o r_0^3}{(r_0+s)^3},
\nonumber
\end{align}
and
\begin{align}
&
\sum_{\substack{m+n=4 \\ 2\geq m'+n \geq 1}} \leVe P_{m,m'+2}( \dcircnabla, \ufls ) \ast \circnabla^{n} \partial_{\us}^{m'} ( \Omega^2 \slashg^{-1} ) \riVe^{n,p}
\leq
\frac{c(n,p) c_o^2 \udelta_o^2 r_0^2}{(r_0+s)^2},
\nonumber
\\
&
\sum_{\substack{m_1+m_2 \leq 4 \\ m_1 \leq 2, m_2 \leq 2}} \leVe P_{m_1,1} ( \dcircnabla, \ufls) \ast P_{m_2,1}(\circnabla, \ufls) \ast (\Omega^2 \slashg^{-1} ) \riVe^{n,p}
\leq
\frac{c(n,p) c_o^2 \udelta_o^2 r_0^2}{(r_0+s)^2}.
\nonumber
\end{align}
The last two estimates follow from the inequality
\begin{align}
\label{eqn 3.32}
\leVe \prod_{i=1}^m h_i \riVe^{n,p} \leq c(m,n,p) \prod_{i=1}^m \leVe h_i \riVe^{n,p},
\quad
n\geq 1,p>2 \text{ or } n\geq 2, p>1.
\end{align}
For example we demonstrate the estimate for the term $2\Omega^2 {\slashg^{-1}}^{ij} \dcircnabla_k \ufls_i \dcircnabla^k \ufls_j$. By bootstrap assumption \ref{assum 3.4},
\begin{align}
\leVe \dcircnabla_k \ufls_i  \riVe^{n,p} \leq c_o \udelta_o r_o.
\nonumber
\end{align}
Further by the bounds of $\Omega$, $\slashg$ and their derivatives, we have
\begin{align}
\leVe \Omega^2 {\slashg^{-1}}^{ij}  \riVe^{n,p} \leq \frac{c(n,p)}{(r_0 + s)^2}.
\nonumber
\end{align}
Thus by inequality \eqref{eqn 3.32}, we obtain
\begin{align}
\leVe 2\Omega^2 {\slashg^{-1}}^{ij} \dcircnabla_k \ufls_i \dcircnabla^k \ufls_j \riVe^{n,p}
\leq
\frac{c(n,p) c_o^2 \udelta_o^2 r_0^2}{(r_0+s)^2}.
\nonumber
\end{align}
Therefore we prove estimates \eqref{eqn 3.23} \eqref{eqn 3.24} for $\Xls$ and $\rels$. Now we can use these estimates to improve the estimates for $\dslashd \ufls$ and $\overline{\ufls}^{\circg}$ in the bootstrap assumption.

For $\dslashd \ufls$, we apply Gronwall's inequality lemma \ref{lem 3.5} stated later to equation \eqref{eqn 3.20}. Assume that $\delta$ is sufficiently small such that
\begin{align}
\leVe \Xls \riVe^{n+1,p}
\leq
&
\frac{c(n,p)( c_{m,m} \eps \udelta_m + c_{m,o} \eps \udelta_o^2 + c_o \eps \udelta_o )r_0^2}{(r_0+s)^3}
+
\frac{c(n,p) c_o \udelta_o r_0}{(r_0+s)^2}
\nonumber
\\
\leq
&
\left\{
\frac{c(n,p)( c_{m,m} \delta + c_{m,o} \delta^3 + c_o \delta^2 )r_0}{r_0+s}
+
c(n,p) c_o \delta
\right\}
\frac{r_0}{(r_0+s)^2}
\nonumber
\\
\leq
&
\frac{ r_0}{(r_0+s)^2},
\nonumber
\end{align}
then by lemma \ref{lem 3.5} we obtain
\begin{align}
\leVe \dcircDelta \ufls \riVe^{n,p} 
&
\leq c(n,p) \left\{ \leVe \dcircDelta \ufl{s=0} \riVe^{n,p} + \int_0^s \leVe \rel{s=t} \riVe^{n,p} \d t \right\}
\nonumber
\\
&
\leq 
c(n,p) \udelta_o r_0
+
c(n,p)( c_{m,m} \eps \udelta_m + c_{m,o} \eps \udelta_o^2 + c_o \eps \udelta_o )c_o \udelta_o r_0
+
c(n,p) c_o^2 \udelta_o^2 r_0
\nonumber
\\
&
\leq
c(n,p) \udelta_o r_0
+
c(n,p) ( c_{m,m} \delta + c_{m,o} \delta^3 + c_o \delta^2 + c_o \delta) c_o \udelta_o r_0.
\nonumber
\end{align}
Thus by the elliptic estimate on the sphere, we obtain that
\begin{align}
\label{eqn 3.38}
\leVe \dslashd \ufls \riVe^{n+1,p} 
\leq 
c(n,p) \udelta_o r_0
+
c(n,p) ( c_{m,m} \delta  + c_{m,o} \delta^3 + c_o \delta^2 + c_o \delta) c_o \udelta_o r_0.
\end{align}

For $\overline{\ufls}^{\circg}$, by equation \eqref{eqn 3.12} and the bootstrap assumption, we have that
\begin{align}
\leve \dpartial_s \ufls \rive 
\leq &
\leve b^i \dpartial_i \ufls \rive + \leve \Omega^2 \left( \slashg^{-1} \right)^{ij} \ufls_i \ufls_j \rive
\nonumber
\\
\leq &
\frac{c(n,p)( c_{m,m} \eps \udelta_m + c_{m,o} \eps \udelta_o^2 + c_o \eps \udelta_o )c_o \udelta_o r_0^3}{(r_0+s)^3}
+
\frac{c(n,p) c_o^2 \udelta_o^2 r_0^2}{(r_0+s)^2},
\label{eqn 3.39}
\end{align}
therefore we have that
\begin{align}
\leve \overline{\ufls}^{\circg} \rive
\leq & 
\leve \overline{\ufl{s=0}}^{\circg} \rive+
\int_{0}^s \leve \overline{\dpartial_s \ufl{s=t}}^{\circg} \rive \d t
\nonumber
\\
\leq &
\udelta_m r_0
+
c(n,p)( c_{m,m} \eps \udelta_m + c_{m,o} \eps \udelta_o^2 + c_o \eps \udelta_o )c_o \udelta_o r_0
+
c(n,p) c_o^2 \udelta_o^2 r_0
\nonumber
\\
\leq &
\left[ 1+ c(n,p) c_{m,m} c_o \delta^2 \right] \udelta_m r_0
+
c(n,p) \left( c_{m,o} c_o \delta^2 + c_o^2 \delta + c_o^2 \right) \udelta_o^2 r_0.
\label{eqn 3.40}
\end{align}

Now we choose the constants $c_o, c_{m,m}, c_{m,o}$ and $\delta$ such that \eqref{eqn 3.38} \eqref{eqn 3.40} improve the estimates for $\ufls$ in the bootstrap assumption. First in the process to obtain estimates \eqref{eqn 3.38} \eqref{eqn 3.40}, we assumed that
\begin{align}
&
\delta \leq \frac{1}{2}, 
\quad
(c_o+c_{m,m} +c_{m,o}) \delta \leq 1,
\nonumber
\\
&
\frac{c(n,p)(c_{m,m} \delta  + c_{m,o} \delta^3 + c_o \delta^2 )r_0}{r_0+s}
+
c(n,p) c_o \delta
\leq
1.
\nonumber
\end{align}
Second, we impose the following inequalities which are sufficient for improving the bootstrap assumption by estimates \eqref{eqn 3.38} \eqref{eqn 3.40}
\begin{align}
c(n,p) 
+
c(n,p) (c_{m,m} \delta  + c_{m,o} \delta^3 + c_o \delta^2 + c_o \delta) c_o
<
c_o,
\nonumber
\\
1+ c(n,p) c_{m,m} c_o \delta^2 < c_{m,m},
\nonumber
\\
c(n,p) \left(  c_{m,o} c_o \delta^2 + c_o^2 \delta + c_o^2 \right) < c_{m,o}.
\nonumber
\end{align}
There are positive solutions for $c_o, c_{m,m}, c_{m,o}$ and $\delta$ to the above inequalities, for example we can first choose
\begin{align}
c_o =2 c(n,p) , \quad c_{m,m} =2, \quad c_{m,o} = 8 c(n,p)^3,
\nonumber
\end{align}
then choose $\delta$ sufficiently small such that the above inequalities hold.

Therefore we show that there exist constants $c_o, c_{m,m}, c_{m,o}$ and $\delta$ depending on $n,p$, such that if bootstrap assumption \ref{assum 3.4} hold, we can improve the estimates of $\ufls$ to extend the assumption further beyond the interval. Thus by the bootstrap argument, the theorem is proved for these constants $c_o, c_{m,m}, c_{m,o}$ and $\delta$ depending on $n,p$.
\end{proof}

For the rest of this paper, we will simply write $c_o, c_{m,m}, c_{m,o}$ as $c$ or $c(n,p)$. In the above proof of theorem \ref{thm 3.3}, we apply the following Gronwall's inequality, which is proved in the appendix.
\begin{lemma}[Gronwall's inequality]\label{lem 3.5}
Consider the cylinder $\mathbb{S}^2 \times \mathbb{R}$. Let $s$ be the parameter in $\mathbb{R}$ and $\vartheta$ or $(\theta^1, \theta^2)$be the variables on the sphere. Suppose that $\mathbb{S}^2$ is endowed with the round metric $\circg$. Let $\Xls$ be a $\mathbb{S}^2$ tangential vector field and $\rels$ be a function on the cylinder. Assume that
\begin{align}
\Xls \in \mathrm{W}^{n+1,p}\left( \mathbb{S}^2 \right),
\nonumber
\end{align}
and $\Xls$ satisfies the estimate
\begin{align}
\leVe \Xls \riVe^{n+1,p} \leq \frac{k r_0}{(r_0+s)^2},
\nonumber
\end{align}
where $n\geq 1, p>2$ or $n\geq 2, p > 1$.
Then there exists a constant $c(n,p)$ such that for any solution $\ufls$ of the following equation
\begin{align}
\partial_s \ufls + \Xls^i \partial_i \ufls = \rels,
\nonumber
\end{align}
we have that
\begin{align}
\leVe \ufls \riVe^{m,p} \leq \exp (c(n,p) k) \left\{ \leVe \ufl{s=0} \riVe^{m,p} + \int_0^s \leVe \rel{s=t} \riVe^{m,p} \d t \right\},
\nonumber
\end{align}
for all integers $0\leq m\leq n+1$.
\end{lemma}

\begin{corollary}\label{coro 3.6}
Under the assumptions of theorem \ref{thm 3.3},
there exists a small positive number $\delta$ depending on $n,p$, that if $\eps,\udelta_o,\udelta_m$ are suitably bounded such that $\eps, \udelta_o, \eps\udelta_m\leq \delta$, then $\ufls$ satisfies the following estimates
\begin{align}
&
\leve \overline{\ufls}^{\circg} - \overline{\ufl{s=0}}^{\circg} \rive 
\leq
c(n,p) (\eps \udelta_m \udelta_o + \udelta_o^2) r_0,
\label{eqn 3.19}
\end{align}
for $s\in(s_{\alpha}, s_{\beta})$.

Let $a$ be a positive number less than $1$. There is a constant $c$ depending on $n,p$, that if $\udelta_m + c \udelta_o < a \kappa$, then $\leve \ufls \rive<a \kappa r_0$ for all $s\in (-\kappa r_0, +\infty)$.
\end{corollary}
\begin{proof}
Estimate \eqref{eqn 3.19}
follows from estimate \eqref{eqn 3.39}
for $\dpartial_s \ufls$ and
\begin{align*}
\leve \overline{\ufls}^{\circg} - \overline{\ufl{s=0}}^{\circg} \rive \leq  \int_0^s \leve \overline{\dpartial_s \ufl{s=t}}^{\circg} \rive \d t.
\end{align*}

Then from estimates \eqref{eqn 3.14} \eqref{eqn 3.19}, we obtain that
\begin{align*}
\leve \ufls \rive 
\leq  
&
\leve \overline{\ufl{s=0}}^{\circg} \rive
+
\leve \overline{\ufls}^{\circg} - \overline{\ufl{s=0}}^{\circg} \rive 
+
c(n,p) \leVe \dslashd \ufls \riVe^{n+1,p}
\\
\leq
&
\udelta_m r_0 + c(n,p) (\eps \udelta_m \udelta_o + \udelta_o^2) r_0 + c(n,p) c_o \udelta_o r_0
\\
\leq
&
\udelta_m r_0 + ( 2 c(n,p) \delta +c(n,p) c_o ) \udelta_o r_0
\end{align*}
Therefore we choose $c=2 c(n,p) \delta +c(n,p) c_o$, then if $\udelta_m + c \udelta_o<a \kappa$, we have $\leve \ufls \rive < a \kappa r_0$. Thus we prove that given $\udelta_m + c \udelta_o<a \kappa$, if $\leve \ufls \rive < \kappa r_0$ on $(s_{\alpha}, s_{\beta})$, then $\leve \ufls \rive < a \kappa r_0$ on $(s_{\alpha}, s_{\beta})$. Hence by a bootstrap argument, we show that if $\udelta_m + c \udelta_o<a \kappa$, then $\leve \ufls \rive < a \kappa r_0$.
\end{proof}

\begin{remark}
In the proof of theorem \ref{thm 3.3}, we use the $L^{\infty}$ bounds of the metric components $\Omega^2$, $\vec{b}$, $\slashg$ and their derivatives up to the order $n+2$. No bounds on the derivatives in the direction $\partial_s$ are required for theorem \ref{thm 3.3}.
\end{remark}

\section{Perturbations of incoming null hypersurfaces}\label{sec 4}
\noindent
In this section, we introduce the perturbation function $\dd{\ufls}$ to describe the perturbation of null hypersurfaces. Then we estimate $\dd{\ufls}$ when given the initial perturbation function $\dd{\ufl{s=0}}$.

\subsection{Perturbation functions}\label{subsec 4.1}
Suppose that $\ucalH_1$ and $\ucalH_2$ are two incoming null hypersurfaces in $(M,g)$. Let $\{\Sigmal{a}_s\}$ be the corresponding foliation on $\ucalH_a, a=1,2$. Assume that $\Sigmal{a}_s$ is parametrised by $\left( s, \ufl{a,s} \right)$. We define the perturbation function of the parametrisation functions $\ufl{a,s}$ as
\begin{align}
\dd{\ufls} = \ufl{2,s} -\ufl{1,s}.
\nonumber
\end{align}

From equations \eqref{eqn 3.12} \eqref{eqn 3.16} \eqref{eqn 3.20} on the propagation of the parametrisation function, we derive the following equations for the perturbation function $\dd{\ufls}$
\begin{align}
\left.
\begin{aligned}
&
\dpartial_s \dd{\ufls} 
=
\dd{F(s, \ufls, \dslashd \ufls)}
=
F(s, \ufl{2,s}, \dslashd \ufl{2,s}) - F(s, \ufl{1,s}, \dslashd \ufl{1,s} ),
\\
&
F(s, \ufls, \dslashd \ufls)
=
-b^i \dpartial_i \ufls + \Omega^2 \left( \slashg^{-1} \right)^{ij} \dpartial_i \ufls \dpartial_j \ufls,
\end{aligned}
\right.
\label{eqn 4.2}
\end{align}
and
\begin{align}
&
\dpartial_s \Big( \dcircDelta \dd{\ufls} \Big)
=
\Xl{1,s}^i \dpartial_i \Big( \dcircDelta \dd{\ufls} \Big) + \dd{\Xls^i} \dpartial_i \Big( \dcircDelta \ufl{2,s} \Big) + \dd{\rels},
\label{eqn 4.3}
\end{align}
where $\dd{\Xls}$, $\dd{\rels}$ are the differences of the corresponding quantities on $\ucalH_a$. In the above equations, we view $\Xls$ and $\rels$ as functions of $s, \ufls, \dslashd \ufls, \dcircnabla^2 \ufls$, i.e.
\begin{align}
\Xl{a,s} = \Xls(s, \ufl{a,s}, \dslashd \ufl{a,s}),
\quad
\rels = \rels(s, \ufl{a,s}, \dslashd \ufl{a,s}, \dcircnabla^2 \ufl{a,s}),
\label{eqn 4.4}
\end{align}
by the formulae of their definitions. We shall use equations \eqref{eqn 4.2} \eqref{eqn 4.3} to obtain the estimates for the perturbation functions $\dd{\ufls}$.
\begin{remark}\label{rem 4.1}
We shall see that the estimates obtained by integrating equation \eqref{eqn 4.3} will lose one order of derivative in the regularity of $\dd{\ufls}$. It is caused by the term $\dd{\Xls^i} \dpartial_i \Big( \dcircDelta \ufl{2,s} \Big)$, which involves the third order derivatives of $\ufl{2,s}$.
\end{remark}

\subsection{Estimates for perturbation functions}
We shall prove the following theorem estimating the perturbation functions $\dd{\ufls}$.
\begin{theorem}\label{thm 4.2}
Let $\ucalH_a, a=1,2$ be two incoming null hypersurfaces and $\{\Sigmal{a}_s\}$ be the corresponding foliation on $\ucalH_a$. Suppose that $\Sigmal{s}_a$ is parametrised by $(s, \ufl{a,s})$. Let $\dd{\ufls}$ be the perturbation function for $\ucalH_a$. Let $(s_{\alpha},s_{\beta})$ be an interval containing $s=0$. Suppose that $\leve \ufl{a,s} \rive <\kappa r_0$ for all $s \in (s_{\alpha},s_{\beta})$.

Assume that the parametrisation functions $\ufl{a,s=0}$ and perturbation function $\dd{\ufl{s=0}}$ satisfy the following estimates
\begin{align}
&
\leVe \dslashd \ufl{s=0} \riVe_{\circg}^{n+1,p} \leq \udelta_o r_0,
\quad
\leve \overline{\ufl{s=0}}^{\circg} \rive \leq  \udelta_m r_0,
\nonumber
\\
&
\leVe \dslashd \dd{\ufl{s=0}} \riVe^{n,p} \leq \ufrakd_o r_0,
\quad
\leve \overline{\dd{\ufl{s=0}}}^{\circg} \rive \leq  \ufrakd_m r_0,
\nonumber
\end{align}
where $n\geq 1, p>2$ or $n\geq 2, p>1$.

There exists a small positive constant $\delta$ depending on $n,p$, that if $\eps, \udelta_o, \udelta_m, \ufrakd_o, \ufrakd_m$ are suitably bounded such that $\eps, \udelta_o, \ufrakd_o, \eps \udelta_o, \eps \ufrakd_m \leq \delta$, then there exist constants $c_o, c_{o,m}, c_m, c_{m,o}$ depending on $n,p$ that $\dd{\ufls}$ satisfies the following estimates
\begin{align}
&
\leVe \dslashd \dd{\ufls} \riVe^{n,p}
\leq
c_o \ufrakd_o r_0 + c_{o,m} ( \udelta_o^2 + \eps \udelta_o ) \ufrakd_m r_0,
\label{eqn 4.7}
\\
&
\leve \overline{\dd{\ufls}}^{\circg} \rive
\leq
c_m \ufrakd_m r_0 + c_{m,o} ( \udelta_o + \eps \udelta_m) \ufrakd_o r_0,
\label{eqn 4.8}
\end{align}
for all $s\in(s_{\alpha},s_{\beta})$.
\end{theorem}

We shall prove theorem \ref{thm 4.2} by bootstrap arguments. Introduce the following bootstrap assumption.
\begin{assumption}[Bootstrap assumption of $\dd{\ufls}$ for theorem \ref{thm 3.3}]\label{assum 4.3}
Estimates \eqref{eqn 4.7} \eqref{eqn 4.8} of $\dd{\ufls}$ hold for all $s\in [0,s_a]$ or $s\in[s_a,0]$.
\end{assumption}

We shall obtain estimates of $\dd{F(s,\ufls, \dslashd \ufls)}$, $\dd{\Xls}$ and $\dd{\rels}$ under bootstrap assumption \ref{assum 4.3}, then use equations \eqref{eqn 4.2}, \eqref{eqn 4.3} and Gronwall's inequality to improve estimates of $\dd{\ufls}$ such that bootstrap assumption \ref{assum 4.3} can be extended to a slightly larger interval by continuity. The constants $c_o, c_{o,m},c_m, c_{m,o}$ are determined in the end of the bootstrap argument.

It is necessary to point out that the bootstrap argument is not really needed for the proof of theorem \ref{thm 4.2}, as the permutation function $\dd{\ufls}$ appears in equations \eqref{eqn 4.2} \eqref{eqn 4.3} linearly. It is in contrast to the proof of theorem \ref{thm 3.3} where the bootstrap argument is indeed needed.

\begin{proof}[Proof of theorem \ref{thm 4.2}]
Like in the proof of theorem \ref{thm 3.3}, we use $c(n,p)$ or simply $c$ to denote constants depending only on $n,p$ which are not necessary the same throughout the proof. We assume that $\eps, \udelta_o, \ufrakd_o, \eps \udelta_o, \eps \ufrakd_m \leq \delta \leq \frac{1}{2}$ and moreover $\delta$ is sufficiently small that
\begin{align}
(c_o + c_{o,m} + c_m + c_{m,o}) \delta \leq 1.
\nonumber
\end{align}

To simplify the notations, let $\ud_o, \ud_m$ denote
\begin{align}
\ud_o = c \udelta_o,
\quad
\ud_m = c (\udelta_m + \udelta_o^2),
\nonumber
\end{align}
obtained in theorem \ref{thm 3.3}. Then by theorem \ref{thm 3.3}
\begin{align}
\leVe \dslashd \ufl{a,s} \riVe^{n+1,p} \leq \ud_o r_0,
\quad
\leve \overline{\ufl{a,s}}^{\circg} \rive \leq \ud_m r_0. 
\nonumber
\end{align}

Let $\ubfd_o, \ubfd_m$ denote
\begin{align}
\ubfd_o = c_o \ufrakd_o + c_{o,m} (\udelta_o^2 + \eps \udelta_o) \ufrakd_m,
\quad
\ubfd_m = c_m \ufrakd_m + c_{m,o} ( \udelta_o + \eps \udelta_m) \ufrakd_o.
\nonumber
\end{align}
in estimates \eqref{eqn 4.7}, \eqref{eqn 4.8}. Then by bootstrap assumption \ref{assum 4.3}, for $s\in[0,s_a]$ or $s\in [s_a,0]$, perturbation function $\dd{\ufls}$ satisfies
\begin{align}
\leVe \dcircnabla^{m+1} \dd{\ufls} \riVe^{n-m,p}
\leq
\ubfd_o r_0,
\quad
\leve \overline{\dd{\ufls}}^{\circg} \rive
\leq
\ubfd_m r_0.
\nonumber
\end{align}

By estimates obtained in the proof of theorem \ref{thm 3.3}, $F(s, \ufl{a,s}, \dslashd \ufl{a,s})$, $\Xl{a,s}$ and $\rel{a,s}$ satisfy
\begin{align}
&
\leve F(s, \ufl{a,s}, \dslashd \ufl{a,s}) \rive
\leq
\frac{c \eps r_0^3}{(r_0+s)^3} \left( \ud_m + \ud_o \right) \ud_o
+
\frac{c r_0^2}{(r_0+s)^2} \ud_o^2,
\label{eqn 4.14}
\\
&
\leVe \Xl{a,s} \riVe^{n+1,p}
\leq
\frac{c \eps r_0^2}{(r_0+s)^3} \left( \ud_m + \ud_o \right)
+
\frac{c r_0}{(r_0+s)^2} \ud_o,
\label{eqn 4.15}
\\
&
\leVe \rel{a,s} \riVe^{n,p}
\leq
\frac{c \eps r_0^3}{(r_0+s)^3} \left( \ud_m + \ud_o \right) \ud_o
+
\frac{c r_0^2}{(r_0+s)^2} \ud_o^2.
\label{eqn 4.16}
\end{align}
By bootstrap assumption \ref{assum 4.3} and the bounds in definition \ref{def 2.2}, we obtain estimates for  their perturbations $\dd{F(s,\ufls, \dslashd \ufls)}$, $\dd{\Xls}$ and $\dd{\rels}$
\begin{align}
&
\leve \dd{F(s,\ufls, \dslashd \ufls)} \rive 
\leq 
\left[ \frac{ c \eps \ud_o r_0^3}{(r_0+s)^3} + \frac{ c \ud_o^2 r_0^2}{ (r_0+s)^2} \right] \ubfd_m
+
\left[ \frac{ c \eps \ud_m r_0^3}{(r_0+s)^3} + \frac{ c \ud_o r_0^2}{ (r_0+s)^2} \right] \ubfd_o,
\label{eqn 4.17}
\\
&
\leVe \dd{\Xls} \riVe^{n,p}
\leq
\left[ \frac{c \eps r_0^3}{(r_0+s)^3} + \frac{ c \ud_o r_0^2}{ (r_0+s)^2} \right] \ubfd_m
+
\frac{ c r_0^2}{ (r_0+s)^2} \ubfd_o,
\label{eqn 4.18}
\\
&
\leVe \dd{\rels} \riVe^{n-1,p}
\leq 
\left[ \frac{ c \eps \ud_o r_0^3}{(r_0+s)^3} + \frac{ c \ud_o^2 r_0^2}{ (r_0+s)^2} \right] \ubfd_m
+
\left[ \frac{ c \eps \ud_m r_0^3}{(r_0+s)^3} + \frac{ c \ud_o r_0^2}{ (r_0+s)^2} \right] \ubfd_o.
\label{eqn 4.19}
\end{align}

Estimates \eqref{eqn 4.17}-\eqref{eqn 4.19} can be obtained from \eqref{eqn 4.14}-\eqref{eqn 4.16} via the following heuristic formal process:
\begin{align}
&
\text{type 1 term}: 
&&
\frac{\eps^{i} r_0^j}{(r_0+s)^k}
&&
\longrightarrow
&&
\frac{\eps^{i} r_0^j}{(r_0+s)^k} (\ubfd_m+ \ubfd_o),
\nonumber
\\
&
\text{type 2 term}: 
&&
\ud_m
&&
\longrightarrow
&&
\ubfd_m,
\nonumber
\\
&
\text{type 3 term}: 
&&
\ud_o
&&
\longrightarrow
&&
\ubfd_o,
\nonumber
\end{align}
and for a term which is a polynomial of the above terms, we can apply the Leibniz rule to it, which is demonstrated by the following diagram.
\begin{center}
\begin{tikzcd}
&
\frac{\eps^i r_0^j}{(r_0+s)^k} \ud_m^l \ud_o^q
\arrow[dl]
\arrow[d]
\arrow[dr]
&
\\
\left( \frac{\eps^i r_0^j}{(r_0+s)^k} (\ubfd_m + \ubfd_o) \right) \ud_m^l \ud_o^q 
&
\frac{\eps^i r_0^j}{(r_0+s)^k} \ud_m^{l-1} \ubfd_m \ud_o^q
&
\frac{\eps^i r_0^j}{(r_0+s)^k} \ud_m^l \ud_o^{q-1} \ubfd_o
\end{tikzcd}
\end{center}
Estimates \eqref{eqn 4.17}-\eqref{eqn 4.19} obey the above pattern. The rigorous proofs follow from the following identity
\begin{align}
\dd{\prod_{i=1}^m h_i} 
=
\dd{h_1} \hle{2}_2 \cdots \hle{2}_m 
+
\hle{1}_1 \dd{h_2} \hle{2}_3 \cdots \hle{2}_m 
+
\cdots
+
\hle{1}_1 \cdots \hle{1}_{m-1} \dd{h_m},
\nonumber
\end{align}
and the following inequality on the Sobolev norm for products of functions
\begin{align}
\leVe \prod_{i=1}^m h_i \riVe^{n-1,p} \leq c(m,n,p) \leVe h_1 \riVe^{n-1,p} \cdot \prod_{i=2}^m \leVe h_i \riVe^{n,p},
\quad
n\geq 1,p>2 \text{ or } n\geq 2, p>1.
\nonumber
\end{align}
For $\dd{\Xls}$, the term with worst regularity is the term involving $\dd{\dslashd \ufls}$ from
$$\dd{\Omega^2 \left( \slashg^{-1} \right)^{ij} \dpartial_j \ufls}.$$
For $\dd{\rels}$, the terms with worst regularity are the terms involving $\dd{\dcircnabla^2_{kl} \ufls}, \dd{\dcircDelta \ufls}$ from
$$
\dd{(\circnabla^k b^i) \dcircnabla^2_{ki} \ufls},
\quad
\dd{\partial_{\us} b^i\ \ufls_i \dcircDelta \ufls},
\quad
\dd{\partial_{\us}b^i\ \ufls^k \dcircnabla^2_{ki} \ufls },
$$
$$
\dd{\circnabla^k( \Omega^2 {\slashg^{-1}})^{ij}\ \ufls_i \dcircnabla^2_{kj} \ufls},
\quad
\dd{\Omega^2 {\slashg^{-1}}^{ij} \dcircnabla_k \ufls_i \dcircnabla^k \ufls_j},
$$
$$
\dd{\partial_{\us}(\Omega^2 \slashg^{-1})^{ij} \ufls_i \ufls^k \dcircnabla^2_{kj} \ufls},
\quad
\dd{\partial_{\us} ( \Omega^2 \slashg^{-1})^{ij} \ufls_i \ufls_j  \dcircDelta \ufls}.
$$

We integrate equations \eqref{eqn 4.2}, \eqref{eqn 4.3} and use the bounds \eqref{eqn 4.17}-\eqref{eqn 4.19} to obtain estimates of $\dd{\ufls}$.

For $\overline{\dd{\ufls}}^{\circg}$, we have
\begin{align}
\leve \overline{\dd{\ufls}}^{\circg} \rive 
\leq
&
\leve \overline{\dd{\ufl{s=0}}}^{\circg} \rive
+
\int_0^s \leve \overline{\dd{F(t,\ufl{s=t}, \dslashd \ufl{s=t})}}^{\circg} \rive \d t
\nonumber
\\
\leq 
&
\ufrakd_m r_0
+
c \left( \eps \ud_o + \ud_o^2  \right) \ubfd_m r_0
+
c \left( \eps \ud_m + \ud_o \right) \ubfd_o r_0.
\label{eqn 4.25}
\end{align}

For $\dslashd \dd{\ufls}$, we apply Gronwall's inequality in \ref{lem 3.5} to obtain
\begin{align}
\leVe \dcircDelta \dd{\ufls} \riVe^{n-1,p}
\leq
&
c \leVe \dcircDelta \dd{\ufl{s=0}} \riVe^{n-1,p}
+
c \int_0^s  \leVe \dd{\Xl{s=t}^i} \dpartial_i \Big( \dcircDelta \ufl{2,s=t} \Big)  \riVe^{n-1,p} \d t
\nonumber
\\
&
+
c \int_0^s \leVe \dd{\rel{s=t}} \riVe^{n-1,p} \d t
\nonumber
\\
\leq
&
c \ud_o r_0
+
c \left( \eps \ud_o + \ud_o^2  \right) \ubfd_m r_0
+
c \left( \eps \ud_m + \ud_o \right) \ubfd_o r_0.
\nonumber
\end{align}
Then by the theory of elliptic equations on the sphere, we obtain that
\begin{align}
\leVe \dslashd \dd{\ufls} \riVe^{n,p}
\leq
&
c\ufrakd_o r_0
+
c \left( \eps \ud_m + \ud_o \right) \ubfd_o r_0
+
c \left( \eps \ud_o + \ud_o^2  \right) \ubfd_m r_0.
\label{eqn 4.27}
\end{align}
Substituting $\ud_m, \ud_o, \ubfd_m, \ubfd_o$ to estimates \eqref{eqn 4.25} \eqref{eqn 4.27}, we obtain
\begin{align}
\leve \overline{\dd{\ufls}}^{\circg} \rive
\leq
&
\ufrakd_m r_0 + c \left( \eps \udelta_o +  \udelta_o^2 \right) \left[ c_m \ufrakd_m + c_{m,o} ( \udelta_o + \eps \udelta_m) \ufrakd_o\right] r_0
\nonumber
\\
&
+
c \left( \eps \udelta_m + \udelta_o \right) \left[ c_o \ufrakd_o + c_{o,m} (\udelta_o^2 + \eps \udelta_o) \ufrakd_m \right] r_0,
\nonumber
\\
\leVe \dslashd \dd{\ufls} \riVe^{n,p}
\leq
&
c \ufrakd_o r_0 + c \left( \eps \udelta_m + \udelta_o \right) \left[ c_o \ufrakd_o + c_{o,m} (\udelta_o^2 + \eps \udelta_o) \ufrakd_m \right] r_0
\nonumber
\\
&
+
c \left( \eps \udelta_o + \udelta_o^2  \right) \left[ c_m \ufrakd_m + c_{m,o} ( \udelta_o + \eps \udelta_m) \ufrakd_o\right] r_0.
\nonumber
\end{align}

We shall choose constants $c_o, c_{o,m}, c_m, c_{m,o}$ and $\delta$ such that estimates \eqref{eqn 4.25} \eqref{eqn 4.27} improve the estimates of $\dd{\ufls}$ in bootstrap assumption \ref{assum 4.3}. First, we require that
\begin{align}
\delta \leq \frac{1}{2},
\quad
(c_o + c_{o,m} + c_m + c_{m,o}) \delta < 1.
\nonumber
\end{align}
Moreover we require that
\begin{align}
1+ c \left( \eps \udelta_m + \udelta_o^2 \right) c_m + c ( \eps \udelta_m + \udelta_o) (\udelta_o^2 + \eps \udelta_o) c_{o,m} < c_m,
\nonumber
\\
c(\eps \udelta_m + \udelta_o^2) c_{m,o} + c  c_o < c_{m,o},
\nonumber
\\
c+ c(\eps \udelta_m + \udelta_o) c_o + c( \eps \udelta_o + \udelta_o^2) (\udelta_o + \eps \udelta_m) c_{m,o} < c_o,
\nonumber
\\
c( \eps \udelta_m + \udelta_o) c_{o,m} + c c_m < c_{o,m}.
\nonumber
\end{align}
The above system has nonempty solution set. For example, we can choose
\begin{align}
( c_m, c_{m,o}, c_o, c_{o,m} ) = ( 2, c + 2 c^2, 2c, 3 c),
\nonumber
\end{align}
and let $\delta$ sufficiently small such that
\begin{align}
c \left( \eps \udelta_m + \udelta_o^2 \right) c_m + c ( \eps \udelta_m + \udelta_o) (\udelta_o^2 + \eps \udelta_o) c_{o,m} 
<1,
\nonumber
\\
(\eps \udelta_m + \udelta_o^2) c_{m,o} < 1,
\nonumber
\\
(\eps \udelta_m + \udelta_o) c_o + ( \eps \udelta_o + \udelta_o^2) (\udelta_o + \eps \udelta_m) c_{m,o}
<1,
\nonumber
\\
( \eps \udelta_m + \udelta_o) c_{o,m} <1.
\nonumber
\end{align}

Therefore we prove that there exist constants $c_m, c_{m,o}, c_o, c_{o,m}$ and $\delta$ sufficiently small depending on $n,p$ such that we can improve the estimates in bootstrap assumption \ref{assum 4.3} to extend the estimates to a lager interval. Hence by the bootstrap arguments, we prove the theorem for these constants $c_m, c_{m,o}, c_o, c_{o,m}$ and $\delta$.
\end{proof}

In the rest of this paper, we simply denote $c_m, c_{m,o}, c_o, c_{o,m}$ by $c(n,p)$ or just $c$.
\begin{remark}
Apply theorem \ref{thm 4.2}, we can prove the following more precise estimates
\begin{align}
&
\leve \overline{\dd{\ufls} - \dd{\ufl{s=0}}}^{\circg} \rive
\leq
c(\eps \udelta_o + \udelta_o^2) \ufrakd_m r_0
+
c(  \eps \udelta_m + \udelta_o) \ufrakd_o r_0,
\nonumber
\\
&
\leVe \dslashd \left( \dd{\ufls} - \dd{\ufl{s=0}} \right) \riVe^{n,p}
\leq
c(  \eps \udelta_m + \udelta_o) \ufrakd_o r_0
+
c(\eps \udelta_o + \udelta_o^2) \ufrakd_m r_0.
\nonumber
\end{align}
In the above proof of the estimate of $\leVe \dslashd \dd{\ufls} \riVe^{n,p}$, we need the $L^{\infty}$ bounds of the metric components up to the $(n+2)$-th order derivatives, since we need to estimate terms like $\leVe \dd{\circDelta b} \riVe^{n-1,p}$, $\leVe \dd{\partial_{\us} \circnabla b} \riVe^{n-1,p}$, $\leVe \dd{\partial_{\us}^2 b} \riVe^{n-1,p}$ and similar $W^{n-1,p}$ norms of the perturbations of second order derivatives of $\Omega^2 \slashg$ in $\leVe \dd{\rels} \riVe^{n-1,p}$.
\end{remark}

\begin{remark}
Note that the estimate of the perturbation function $\dd{\ufls}$ in theorem \ref{thm 4.2} has one less order of derivative than $\ufls$ in theorem \ref{thm 3.3}, as we mentioned in remark \ref{rem 4.1}. If we assume that $\ufl{2,s}$ is constant, which means that the second null hypersurface $\ucalH_2$ is simply one of $\uC_{\us}$ in the double null foliation, then the term $\dd{\Xls^i} \dpartial_i \left( \dcircDelta \ufl{2,s} \right)$ vanishes, thus no loss of derivative occurs in this case and we have the following improved estimate of $\dd{\ufls}$.
\end{remark}

\begin{theorem}\label{thm 4.6}
Under the assumptions of theorem \ref{thm 4.2}, if we assume the following additional assumptions that
\begin{align}
\leVe \dslashd \dd{\ufl{s=0}} \riVe^{n+1,p} \leq \ufrakd_o r_0,
\quad
\ufl{2,s} \equiv \mathrm{const},
\nonumber
\end{align}
then we can obtain the estimate of $\leVe \dslashd \dd{\ufls} \riVe^{n+1,p}$ which takes the same form as in theorem \ref{thm 4.2}.
\end{theorem}
The proof of theorem \ref{thm 4.6} follows the same route as proving theorem \ref{thm 4.2}. The improvement for the regularity of $\dd{\ufls}$ comes from the better assumption of $\dd{\ufl{s=0}}$ and the vanishing of the term $\dd{\Xls^i} \dpartial_i \left( \dcircDelta \ufl{2,s} \right)$. The improved theorem \ref{thm 4.6} requires the $L^{\infty}$ bounds of the metric components up to the $(n+2)$-th order derivatives, no more than the requirement of theorem \ref{thm 4.2}, since we have $\dd{\rels} = -\rel{1,s}$ which doesnot involve terms like $\dd{\circDelta b}, \dd{\partial_{\us} \circnabla b} ,\dd{\partial_{\us}^2 b}$, the perturbations of second order derivatives of metric components.

\section{Linearised perturbation of parametrisation functions}\label{section 5}
\noindent
In this section, we construct the linearised perturbation for parametrisation functions of incoming null hypersurfaces. We will obtain the error estimates between the linearised perturbation and the actual perturbation.

\subsection{Linearised perturbation functions}\label{subsec 5.1}
We consider the same perturbation of incoming null hypersurfaces $\ucalH_a, a=1,2$ as in section \ref{subsec 4.1}
We also constructed the perturbation function $\dd{\ufls}$ in section \ref{subsec 4.1}.
It satisfies equations \eqref{eqn 4.2} \eqref{eqn 4.3}. To construct a reasonable linearised perturbation function for this perturbation, we shall construct the corresponding linearisation of equations \eqref{eqn 4.2} \eqref{eqn 4.3}, then define the linearised perturbation function as the solution of the linearised equations.

We denote the linearised perturbation function by $\bdd{\ufls}$ and define its value at initial time is
\begin{align}
\bdd{\ufl{s=0}} = \ufl{2,s=0} - \ufl{1,s=0}
\nonumber
\end{align}
which is equal to the initial perturbation function $\dd{\ufl{s=0}}$. We construct the following linear equations from equations \eqref{eqn 4.2} \eqref{eqn 4.3}
\begin{align}
&
\dpartial_s \overline{\bdd{\ufls}}^{\circg} = 0
\quad
\Leftrightarrow
\quad
\overline{\bdd{\ufls}}^{\circg} = \overline{\bdd{\ufl{s=0}}}^{\circg},
\label{eqn 5.2}
\\
&
\dpartial_s \Big( \circDelta \bdd{\ufls} \Big) = \Xl{1,s}^i \dpartial_i \Big( \dcircDelta \bdd{\ufls} \Big) - \overline{\Xl{1,s}^i \dpartial_i \Big( \dcircDelta \bdd{\ufls} \Big)}^{\circg}.
\label{eqn 5.3}
\end{align}
We neglect the terms
\begin{align}
\dd{F(s,\ufls, \dslashd \ufls)}, 
\quad
\dd{\Xls^i} \dpartial_i \Big( \dcircDelta \ufl{2,s} \Big),
\quad
\dd{\rels}
\nonumber
\end{align}
in equations \eqref{eqn 4.2} \eqref{eqn 4.3} to obtain the above linear equations. The reason is that these terms are all of quadratic or higher order smallness. For example, from estimates \eqref{eqn 3.39} \eqref{eqn 4.14} \eqref{eqn 4.17}, the terms $F(s,\ufls,\dslashd \ufls)$ and $\dd{F(s, \ufls, \dslashd \ufls)}$ satisfy the estimates
\begin{align}
&
\leve F(s, \ufls, \dslashd \ufls) \rive 
\leq
\frac{c\delta^2 r_0}{(r_0+ s)^2},
\nonumber
\\
&
\leve \dd{F(s,\ufls, \dslashd \ufls)} \rive
\leq
\frac{c \delta r_0^2}{(r_0 + s)^2} \ufrakd_o
+
\frac{c \delta^2 r_0^2}{(r_0 + s)^2} \ufrakd_m,
\nonumber
\end{align}
thus neglecting these terms only cause an error much less than the size of the actual perturbation. We shall see the error between this linearised perturbation function $\bdd{\ufls}$ and $\dd{\ufls}$ in the next section.

Note that a straightforward way to construct a linearised perturbation function is to consider the variation through solutions of equation \eqref{eqn 3.12}. Such a variation, denoted by $v\{\ufls\}$, satisfies the equation of variations corresponding to equation \eqref{eqn 3.12}.
With $\ufl{1,s}$ being the background solution, this equation is
\begin{align*}
\dpartial_{\us} v\{\ufls\}=
& 
- b^i |_{\us=\ufl{1,s}}\cdot \dpartial_i \ \big( v\{\ufls\}  \big)
- \partial_{\us} b^i |_{\us=\ufl{1,s}} \cdot v\{\ufls\}\cdot \dpartial_i \ufl{1,s}
\\
&
+
2\big(\Omega^2 (\slashg^{-1})^{ij} \big)|_{\us=\ufl{1,s}}\cdot \dpartial_i \ufl{1,s}\cdot \dpartial_j \big( v\{\ufls\} \big)
\\
&
+
\partial_{\us} \big(\Omega^2 (\slashg^{-1})^{ij} \big)|_{\us=\ufl{1,s}} \cdot v\{\ufls \} \cdot \dpartial_i \ufl{1,s}\cdot \dpartial_j \ufl{1,s}.
\end{align*}
We see that $\bdd{\ufls}$ and $v\{\ufls\}$ are different, thus $\bdd{\ufls}$ does not represent a variation through solutions of equation \eqref{eqn 3.12}.

\subsection{Estimates for the error of the linearised perturbation}\label{subsec 5.2}
We denote the error of the linearised perturbation function by $\er{\ufls}$ which is the difference of $\dd{\ufls}$ and $\bdd{\ufls}$
\begin{align}
\er{\ufls} = \dd{\ufls} - \bdd{\ufls}.
\nonumber
\end{align}

From equations \eqref{eqn 4.2} \eqref{eqn 4.3} for $\dd{\ufls}$ and equations \eqref{eqn 5.2} \eqref{eqn 5.3}, we derive the following equation for $\er{\ufls}$
\begin{align}
&
\dpartial_s \overline{\er{\ufls}}^{\circg} = \overline{\dd{F(s,\ufls, \dslashd \ufls)}}^{\circg},
\label{eqn 5.8}
\\
&
\dpartial_s \Big( \dcircDelta \er{\ufls} \Big)
=
\Xl{1,s}^i \dpartial_i \Big( \dcircDelta \er{\ufls} \big) + \dd{\Xls^i} \dpartial_i \Big( \dcircDelta \ufl{2,s} \Big) + \dd{\rels}
\nonumber
\\ 
&\hspace{90pt}
+ \overline{\Xl{1,s}^i \dpartial_i \Big( \dcircDelta \bdd{\ufls} \Big)}^{\circg}.
\label{eqn 5.9}
\end{align}
To obtain estimates for $\er{\ufls}$, we still lack the estimate of $\overline{\Xl{1,s}^i \dpartial_i \Big( \dcircDelta \bdd{\ufls} \Big)}^{\circg}$. Its estimate is given by the following lemma.
\begin{lemma}\label{lem 5.1}
Under the assumptions of theorem \ref{thm 4.2}. There exists a constant $c$ and a sufficiently small $\delta$ depending on $n,p$, such that if $\eps, \eps \udelta_m , \udelta_o \leq \delta$, then we have the following estimate for $\dcircDelta{\bdd{\ufls}}$
\begin{align}
\leVe \dcircDelta{\bdd{\ufls}} \riVe^{m,p} \leq c \leVe \dcircDelta{\bdd{\ufl{s=0}}} \riVe^{m,p} 
\nonumber
\end{align}
for all integers $0\leq m \leq n$.
\end{lemma}
\begin{proof}
From estimate \eqref{eqn 4.15} for $\Xl{1,s}$, we have for $\delta$ sufficiently small depending on $n,p$
\begin{align}
\leVe \Xl{1,s} \riVe^{n+1,p} \leq \frac{c \delta r_0}{(r_0+s)^2}.
\nonumber
\end{align}
Then by equation \eqref{eqn 5.9}
and Gronwall's inequality in lemma \ref{lem 3.5}, we have for $\delta$ sufficiently small and any integer $0\leq m \leq n$
\begin{align}
\leVe \dcircDelta{\bdd{\ufls}} \riVe^{m,p} 
\leq
c \left\{ \leVe \dcircDelta{\bdd{\ufl{s=0}}} \riVe^{m,p}  + \int_0^s \Big\vert \overline{\Xl{1,s=t}^i \partial_i \Big( \dcircDelta \bdd{\ufl{s=t}} \Big)}^{\circg} \Big\vert \d t \right\}
\nonumber
\end{align}
For the second term on the right, we have
\begin{align}
\Big\vert \overline{\Xl{1,s}^i \partial_i \Big( \dcircDelta \bdd{\ufl{s=t}} \Big)}^{\circg} \Big\vert
=
&
\Big\vert \overline{ \Big( \dcircdiv \Xl{1,s=t} \Big) \Big( \dcircDelta \bdd{\ufl{s=t}} \Big)}^{\circg} \Big\vert
\nonumber
\\
\leq
&
\leVe \dcircdiv \Xl{1,s=t} \riVe_{L^{\infty}} \leVe \dcircDelta \bdd{\ufl{s=t}} \riVe^{m,p}
\leq
\frac{c\delta r_0}{(r_0+t)^2} \leVe \dcircDelta \bdd{\ufl{s=t}} \riVe^{m,p},
\nonumber
\end{align}
Thus for $\delta$ sufficiently small, we obtain
\begin{align}
\leVe \dcircDelta{\bdd{\ufls}} \riVe^{m,p} \leq c \leVe \dcircDelta{\bdd{\ufl{s=0}}} \riVe^{m,p} .
\nonumber
\end{align}
\end{proof}

Now we can state and prove the following theorem on the estimates for the error of the linearised perturbation function.
\begin{theorem}\label{thm 5.2}
Under the assumptions of theorem \ref{thm 4.2}, there exists a small positive constant $\delta$ depending on $n,p$, that if $\eps$, $\udelta_o$, $\udelta_m$, $\ufrakd_o$, $\ufrakd_m$ are suitably bounded such that $\eps$, $\udelta_o$, $\ufrakd_o$, $\eps \udelta_m$, $\eps \ufrakd_m \leq \delta$, then there exists a constant $c$ depending on $n,p$ such that the error $\er{\ufls}$ of the linearised perturbation satisfies the following estimates
\begin{align}
&
\leVe \dslashd \er{\ufls} \riVe^{n,p} \leq c(\udelta_o + \eps \udelta_m ) \ufrakd_o r_0 + c( \udelta_o^2 + \eps \udelta_o) \ufrakd_m r_0,
\nonumber
\\
&
\leve \overline{\er{\ufls}}^{\circg} \rive 
\leq
c(\udelta_o + \eps \udelta_m ) \ufrakd_o r_0 + c( \udelta_o^2 + \eps \udelta_o) \ufrakd_m r_0,
\nonumber
\end{align}
for $s\in(s_{\alpha}, s_{\beta})$.

Moreover, if we assume the additional assumptions in theorem \ref{thm 4.6},
\begin{align}
\leVe \dslashd \dd{\ufl{s=0}} \riVe^{n+1,p} \leq \ufrakd_o r_0,
\quad
\ufl{2,s} \equiv \mathrm{const},
\nonumber
\end{align}
then we can improve to obtain the estimate of $\leVe \dslashd \er{\ufls} \riVe^{n+1,p}$, which takes the same form as above.
\end{theorem}
\begin{proof}
Since we define the linearised perturbation function $\bdd{\ufl{s=0}}= \ufl{2,s=0} - \ufl{1,s=0}= \dd{\ufl{s=0}}$ at $s=0$, we have
\begin{align}
\er{\ufl{s=0}} =0.
\nonumber
\end{align}
For $\overline{\er{\ufls}}^{\circg}$, we integrate equation \eqref{eqn 5.8} to obtain the estimate. From estimate \eqref{eqn 4.17},
\begin{align}
\leve \dd{F(s, \ufls, \dslashd \ufls)} \rive
\leq
\left[ \frac{ c\eps \udelta_o r_0^3}{(r_0+s)^3} + \frac{c \udelta_o^2 r_0^2}{(r_0+s)^2} \right] \ufrakd_m
+
\left[ \frac{c\eps \udelta_m r_0^3}{(r_0+s)^3} + \frac{c \udelta_o r_0^2}{(r_0+s)^2} \right] \ufrakd_o,
\nonumber
\end{align}
thus we obtain the estimate of $\overline{\er{\ufls}}^{\circg}$
\begin{align}
\leve \overline{\er{\ufls}}^{\circg} \rive 
\leq 
c \left( \eps \udelta_o + \udelta_o^2 \right) \ufrakd_m r_0
+
c \left( \eps \udelta_m  + \udelta_o  \right) \ufrakd_o r_0.
\nonumber
\end{align}
For $\dslashd \er{\ufls}$, we integrate equation \eqref{eqn 5.9} to obtain its estimate. From estimates \eqref{eqn 4.18} \eqref{eqn 4.19},
\begin{align}
&
\leVe \dd{\Xls}^i \dpartial_i \left( \dcircDelta \ufl{2,s} \right) + \dd{\rels} \riVe^{n-1,p}
\nonumber
\\
&
\leq
\left[ \frac{ c\eps \udelta_o r_0^3}{(r_0+s)^3} + \frac{c \udelta_o^2 r_0^2}{(r_0+s)^2} \right] \ufrakd_m
+
\left[ \frac{c\eps \udelta_m r_0^3}{(r_0+s)^3} + \frac{c \udelta_o r_0^2}{(r_0+s)^2} \right] \ufrakd_o.
\nonumber
\end{align}
By lemma \ref{lem 5.1},
\begin{align}
\leve \overline{\Xl{1,s}^i \dpartial_i \left( \dcircDelta \bdd{\ufls} \right)}^{\circg} \rive
=
\leve \overline{\left( \dcircdiv \Xl{1,s} \right) \left( \dcircDelta \bdd{\ufls} \right)}^{\circg} \rive
\leq
\left[ \frac{c\eps\udelta_m  r_0^3}{(r_0+s)^3} + \frac{c\udelta_o r_0^2}{(r_0+s)^2} \right] \ufrakd_o.
\nonumber
\end{align}
Hence by Gronwall's inequality in lemma \ref{lem 3.5}, we obtain that
\begin{align}
\leVe \dcircDelta \er{\ufls} \riVe^{n-1,p}
\leq
c(\eps \udelta_o + \udelta_o^2) \ufrakd_m r_0 + c( \eps \udelta_m + \udelta_o) \ufrakd_o r_0.
\nonumber
\end{align}
Then by the theory of elliptic equations on the shpere, we get that
\begin{align}
\leVe \dslashd \er{\ufls} \riVe^{n,p}
\leq
c(\eps \udelta_o + \udelta_o^2) \ufrakd_m r_0 + c( \eps \udelta_m + \udelta_o) \ufrakd_o r_0.
\nonumber
\end{align}

Under the additional assumptions, by theorem \ref{thm 4.6}, we have improved estimate for $\dslashd \dd{\ufls}$, therefore by the same argument above, we obtain the improved estimate for $\dslashd \er{\ufls}$.
\end{proof}

\section*{Acknowledgements}
\noindent
This paper relaxes the condition of the global existence result of null hypersurfaces in a perturbed Schwarzschild black hole exterior obtained in the author's thesis \cite{L1}. The author is grateful to Demetrios Christodoulou for his generous guidance. The author also thanks Alessandro Carlotto and Lydia Bieri for helpful discussions.

\begin{appendix}

\section{Proof of lemma \ref{lem 3.5}: Gronwall's inequality}
\noindent
We consider the solution of the following propagation equation
\begin{align}
\partial_s \ufls + \Xls^i \partial_i \ufls = \rels
\label{eqn A.1}
\end{align}
and obtain estimates for $\ufls$ by integrating the equation. We first derive the propagation equation for the rotational derivatives of $\ufls$. Let $R_i, i=1,2,3$ be the rotational vector fields on the unit round sphere $(\mathbb{S}^2, \circg)$. In the 3-dimensional Euclidean space, 
$$R_i = \sum_{j,k}\eps_{ijk} x^j \partial_k.$$
Differentiating equation \eqref{eqn A.1}, we get
\begin{align}
\partial_s R_i \ufls + \Xls \left( R_i \ufls \right) + \left[ R_i, \Xls \right] \ufls = R_i \rels.
\nonumber
\end{align}
For the higher order derivatives, we have
\begin{align}
\partial_s \left( R_{i_1} \cdots R_{i_m} \ufls \right) + \Xls \left( R_{i_1} \cdots R_{i_m} \ufls \right)
=
\rel{R_{i_1}, \cdots, R_{i_m},s},
\label{eqn A.3}
\end{align}
where
\begin{align}
\rel{R_{i_1}, \cdots, R_{i_m},s} 
=
&
R_{i_1} \cdots R_{i_n} \rels,
\nonumber
\\
&
-
\sum_{\substack{ 
\{l_1,\cdots, l_s\}\cup \{k_1,\cdots ,k_{n-s}\} \\
=\{1,\cdots, m\}, \\ 
l_1 <\cdots < l_s, \\ 
k_1 < \cdots < k_{n-s}, \\ 
s\leq m-1
}} 
[ R_{i_{k_1}}, \cdots, R_{i_{k_{n-s}}} , \Xls] \left( R_{i_{l_1}} \cdots R_{i_{l_s}} \ufls \right) 
\nonumber
\end{align}
and
\begin{align}
[ R_{i_{k_1}}, \cdots, R_{i_{k_{n-s}}} , \Xls] = \lie_{R_{i_{k_1}}}\cdots \lie_{R_{i_{k_{n-s}}}}\Xls.
\nonumber
\end{align}
The proof of lemma \ref{lem 3.5} relies on the following lemma on the diffeomorphisms generated by the vector field $\Xls$. 
\begin{lemma}\label{lem A.1}
Define the one parameter family of diffeomorphisms $\left\{ \varphi_s \right\}$ generated by $\Xls$, i.e.
\begin{align}
\varphi_s: \mathbb{S}^2 \rightarrow \mathbb{S}^2,
\quad
\frac{\d}{\d s} \varphi_s(\vartheta) = \Xls(\varphi_s(\vartheta)), 
\quad 
\varphi_{s=0}=\mathrm{Id}
\nonumber
\end{align}
We define the push forward metric $g_s$ of $\circg$ via $\varphi_s$
\begin{align}
g_s=(\varphi_s)_* \circg, \quad \dvol_{g_s} = (\varphi_s)_* \dvol_{\circg}.
\nonumber
\end{align}
Then there exists a family of functions $\left\{ \phi_s \right\}$ such that
\begin{align}
\dvol_{g_s} = \phi_s \dvol_{\circg},
\nonumber
\end{align}
and $\phi_s$ satisfies the equation
\begin{align}
\partial_s \phi_s + \Xls^i \partial_i \phi_s = - \phi_s \circdiv \Xls,  \quad  \left( \partial_s  + \Xls^i \partial_i \right) \log \phi_s = - \circdiv \Xls.
\label{eqn A.8}
\end{align}
Assume that 
\begin{align}
\sup_{\mathbb{S}^2} \leve \circdiv \Xls  \rive \leq \frac{k r_0}{(r_0+t)^2},
\nonumber
\end{align}
then
\begin{align}
\leve \log \phi_s \rive \leq k .
\nonumber
\end{align}
\end{lemma}
\begin{proof}
It is sufficient to prove equation \eqref{eqn A.8}. It follows from
\begin{align}
\lie_{\Xls} \dvol_{g_s} = \lie_{\Xls} \left( \phi_s \dvol_{\circg} \right) =0.
\nonumber
\end{align}
\end{proof}

Applying the above lemma, we prove the $L^p$ estimate for $\ufls$.
\begin{lemma}\label{lem A.2}
Under the assumptions of lemma \ref{lem 3.5}, there exists a constant $c$ depending on $p$ such that
\begin{align}
\leVe \ufls \riVe_{L^p}
\leq
\exp (c k ) \left\{ \leVe \ufl{s=0} \riVe_{L^p} + \int_0^s \leVe \rel{s=t} \riVe_{L^p} \d t  \right\}.
\nonumber
\end{align}
\end{lemma}
\begin{proof}
We use the one parameter family defined in lemma \ref{lem A.1}. Denote the pullbacks of $\ufls$ and  $\rels$ via $\varphi_s$ by $\ufls_{\ast}, \rels_{\ast}$,
\begin{align}
\ufls_{\ast} \left(\vartheta \right) = \varphi_s^* (\ufls) (\vartheta) = \ufls \circ \varphi_s (\vartheta),
\quad
\rels_{\ast} \left(\vartheta \right) = \varphi_s^* (\rels) (\vartheta) = \rels \circ \varphi_s (\vartheta).
\nonumber
\end{align}
Then we have the propagation equation
\begin{align}
\frac{\d}{\d s} \ufls_{\ast}(\vartheta) = \rels_{\ast} (\vartheta).
\nonumber
\end{align}
Therefore for $\leVe \ufls_{\ast} \riVe_{L^p}$, we have
\begin{align}
\leve \frac{\d}{\d s} \left( \leVe \ufls_{\ast} \riVe_{L^p} \right)^p \rive
=&
\leve p \int_{\mathbb{S}^2} \rels_{\ast} \cdot \left( \ufls_{\ast} \right)^{p-1} \dvol_{\circg} \rive
\leq
p \leVe \rels_{\ast} \riVe_{L^p} \cdot \leVe \ufls_{\ast} \riVe_{L^p}^{p-1}
\nonumber
\end{align}
Therefore we have
\begin{align}
\leVe \ufls_{\ast} \riVe_{L^p} \leq \leVe \ufl{s=0} \riVe_{L^p} + \int_0^s \leVe \rel{s=t}_{\ast} \riVe_{L^p} \d t. 
\nonumber
\end{align}
Between $L^p$ norms $\leVe \ufls_{\ast} \riVe_{L^p}$ and $\leVe \ufls \riVe_{L^p}$, we have that
\begin{align}
\int_{\mathbb{S}^2} \left( \ufls_{\ast} \right)^p \dvol_{\circg}
=&
\int_{\mathbb{S}^2} \left( (\varphi_s)^{\ast} \ufls \right)^p \dvol_{\circg}
=
\int_{\mathbb{S}^2} \left( \ufls \right)^p \dvol_{\left( \varphi_s \right)_{\ast} \circg}
\nonumber
\\
=&
\int_{\mathbb{S}^2}  \left( \ufls \right)^p \dvol_{g_s}
=
\int_{\mathbb{S}^2} \left( \ufls \right)^p \phi_s \dvol_{\circg}.
\nonumber
\end{align}
Then by lemma \ref{lem A.1}, we have
\begin{align}
\exp(-c k) \leVe \ufls \riVe_{L^p}  \leq \leVe \ufls_{\ast} \riVe_{L^p} \leq \exp (c k) \leVe \ufls \riVe_{L^p} 
\nonumber
\end{align}
since we have the $L^{\infty}$ estimate of $\circdiv \Xls$ from $\leVe \Xls \riVe^{n+1,p}$ by the Sobolev embedding. Similarly, for $\leVe \rels_{\ast} \riVe_{L^p}$ and $\leVe \rels \riVe_{L^p}$, they are also comparable with the constant $\exp (c k)$. Hence
\begin{align}
\leVe \ufls \riVe_{L^p} \leq  \exp(c k) \left\{ \leVe \ufl{s=0} \riVe_{L^p} + \int_0^s \leVe \rel{s=t} \riVe_{L^p} \d t \right\}. 
\nonumber
\end{align}
\end{proof}

Therefore lemma \ref{lem 3.5} follows from equation \eqref{eqn A.3} and lemma \ref{lem A.2}.
\end{appendix}

\end{document}